\documentclass[a4paper,12pt]{article}
\usepackage{amsfonts}
\usepackage{amsmath}
\usepackage{amssymb}
\usepackage{graphicx}
\usepackage{epsf}
\usepackage{epsfig}

\usepackage{color}
\usepackage{afterpage}
\usepackage{verbatim}
\usepackage{amssymb}
\usepackage{graphicx}
\usepackage{amsmath}
\usepackage{amssymb}

\usepackage{graphicx}
\usepackage{epsfig}
\usepackage{labelfig}
\usepackage{amsmath}
\usepackage{color}
\usepackage{afterpage}
\usepackage{verbatim}

\usepackage[small]{caption}

\begin{document}

\title{Domain decomposition finite element/finite difference method
  for the conductivity reconstruction in a hyperbolic equation}

\author{L. Beilina  \thanks{
Department of Mathematical Sciences, Chalmers University of Technology and
Gothenburg University, SE-42196 Gothenburg, Sweden, e-mail: \texttt{\
larisa@chalmers.se}}}

\date{}

\maketitle

\graphicspath{   
{FIGURESArxiv/}
}

\begin{abstract}

We present domain decomposition finite element/finite difference
method for the solution of hyperbolic equation. The domain
decomposition is performed such that finite elements and finite
differences are used in different subdomains of the computational
domain: finite difference method is used on the structured part of the
computational domain and finite elements on the unstructured part of
the domain. The main goal of this method is to combine flexibility of
finite element method and efficiency of a finite difference method.

An explicit discretization schemes for both methods are constructed
such that finite element and finite difference schemes coincide on the
common structured overlapping layer between computational
subdomains. Then the resulting scheme can be considered as a pure
finite element scheme which allows avoid instabilities at the interfaces.

We illustrate efficiency of the domain decomposition method on the
reconstruction of the conductivity function in the hyperbolic equation
in three dimensions.

\end{abstract}

\section{Introduction}

With expanding of new computational technologies and needs of industry
 it is vital importance to use efficient computational methods for
simulation of partial differential equations in two and
three-dimensions when computational domains are very large.  Domain
decomposition methods attracted a lot of interest and is a topic of
current research, see, for example, \cite{Chan, Tos} and references therein.

It is typical that computational domains in industrial applications
often are very large and only some part of this domain presents
interest. In such cases a domain decomposition approach can be
attractive when the simple domain discretization can be used in a
large region and more complex and refined domain discretization is
applied in a smaller part of the domain.  In this paper we propose to
use domain decomposition finite element/finite difference approach for
the solution of hyperbolic equation which combines flexibility of the
finite element method (FEM) and efficiency in terms of speed and
memory usage of finite difference method (FDM).  To do that we extend
a hybrid FEM/FDM method which was developed in \cite{BSA} for the case
of acoustic wave equation, to the case of a more general hyperbolic
equation with two unknown parameters.  Similar to our approach in
\cite{BSA}, we decompose the computational domain such that finite
elements and finite differences are used in different subdomains of
this domain: finite difference method is used in a simple geometry and
finite elements - in the subdomain where we want to get more detailed
information about structure of this domain.  Our goal is to get such a
method which will combine flexibility of finite element method and
efficiency of a finite difference method.  It is well known that the
finite element method allows to get small features of the structure of
the domain through the adaptive mesh refinement. However, this method
is quite computationally expensive comparing with the finite
difference method in terms of time and memory usage, see
\cite{B}-\cite{HybEf} for study of efficiency of these methods.

In this work we derive explicit schemes for both methods such that
finite element and finite difference methods coincide on the common
structured overlapping layer between these domains.  Thus, the
resulting scheme can be considered as a pure finite element scheme
which allows to avoid instabilities at the interfaces  \cite{Brenner}. We implement this
method in an efficient way in the software package WavES \cite{waves} in
C++ using PETSc \cite{petsc} and MPI (message passing interface).

We illustrate the efficiency of the proposed method on the solution of
the hyperbolic coefficient inverse problem (CIP) in three dimensions.
The goal of our numerical simulations is to reconstruct the
conductivity function of the hyperbolic equation from single
observations of the backscattered solution of this equation in space
and time.  We note, that the domain decomposition approach in this
case is particularly feasible for implementing of absorbing boundary
conditions \cite{EM}. To solve our CIP we minimize the corresponding
Tikhonov functional and use Lagrangian approach to do that.  This
approach is similar to one applied recently in \cite{B, BCN, HybEf, BJ} for
the solution of different three-dimensional CIPs: we find optimality
conditions which express stationarity of the Lagrangian, involving the
solution of a state and adjoint equations together with an equation
expressing that the gradient of the Lagrangian with respect to the conductivity
vanishes. Then we construct conjugate gradient algorithm and compute
the unknown conductivity function in an iterative process by solving
in every step of this algorithm the state and adjoint hyperbolic
equations and updating in this way the conductivity function.

We tested our inverse algorithm on the reconstruction of conductivity
function which represents small symmetrical inclusions. This problem
can be interpreted as the problem of the reconstruction of the
symmetrical structure of a waveguide and finding  defects in it.
Our computations show that we can accurately reconstruct large
contrast of the conductivity function as well as location of all  small
inclusions using the domain decomposition method presented in  this work.

The paper is organized as follows. In Section \ref{sec:modelhyb} we
present our mathematical model of hyperbolic equation and in Section
\ref{sec:hyb} we describe domain decomposition approach. Energy
estimate for the equation of Section \ref{sec:modelhyb} is derived in
Section \ref{sec:energyerror1}.  In Section \ref{sec:model} we
formulate state and inverse problems for hyperbolic equation.  The
Tikhonov functional to be minimized and the corresponding Lagrangian
are presented in Section \ref{sec:opt}.  In Section \ref{sec:domdec}
we describe the domain decomposition finite element/finite difference
method to solve the minimization problem of Section \ref{sec:opt}.
Finally, in Section \ref{sec:Numer-Simul} we demonstrate efficiency of
the domain decomposition method on the reconstruction of the
conductivity function in three dimensions.

\section{The mathematical model}

\label{sec:modelhyb}

The model problem in the domain decomposition method is the following hyperbolic equation with the first order absorbing boundary conditions \cite{EM}
\begin{equation}\label{modelhyb}
\begin{split}
\frac{1}{c^2}\frac{\partial^2 u}{\partial t^2}  -  \nabla \cdot ( a \nabla  u)   &= g(x,t),~ \mbox{in}~~ \Omega_T, \\
  u(x,0) = f_0(x), ~~~u_t(x,0) &= f_1(x)~ \mbox{in}~~ \Omega,     \\
\partial _{n}  u& =-\partial _{t} u ~\mbox{on}~ \partial \Omega_T,
\\
\end{split}
\end{equation}
Here $\Omega \subset \mathbb{R}^{3}$ is a convex bounded domain with
the boundary $\partial \Omega \in C^{3}$, $x=(x_1, x_2, x_3) \in
\mathbb{R}^{3}$ and $C^{k+\alpha }$ is H\"older space, where $k\geq 0$
and $\alpha \in \left( 0,\,1\right)$, $c(x), a(x)$ are the wave speed
and the conductivity space-dependent functions, respectively.  We defined by $\Omega_T
:= \Omega \times (0,T), \partial \Omega_T := \partial \Omega \times
(0,T), T > 0$.  We assume that
\begin{equation}
g(x,t) \in L_2(\Omega_T),  f_{0}\in H^{1}(\Omega), f_{1}\in L_{2}(\Omega). \label{f1} 
\end{equation}

For our purpose we use modification of the domain decomposition method
developed in \cite{BSA} which was applied for the solution of the
coefficient inverse problem for the acoustic wave equation in \cite{B,
  BJ}. In this work, we use different version of the method used in
\cite{B, BSA, BJ}, when two functions - the wave speed $c(x)$ and the
conductivity function $a(x)$ - are  introduced in the mathematical model of
the hyperbolic equation.  Similarly with \cite{B}-\cite{BJ},
we decompose $\Omega$ into two open subregions, $\Omega_{FEM}$ and
$\Omega_{FDM}$ such that $\Omega = \Omega_{FEM} \cup \Omega_{FDM}$,
and $\Omega_{FEM} \cap \Omega_{FDM} = \emptyset$.  In $\Omega_{FEM}$
we use finite elements and this domain is such that $\partial
\Omega_{FEM} \subset \partial \Omega_{FDM}$, see figure  
\ref{fig:fig1}.  In $\Omega_{FDM}$ we will use finite difference
method.

We assume that our
coefficients $a(x), c \left( x\right)$ of problem (\ref{modelhyb})   are such
that
\begin{equation} \label{coefic}
\begin{split}
a \left( x\right) &\in \left[ 1,d_1 \right],~~ d_1 = const.>1,~ a(x) =1
\text{ for }x\in  \Omega _{FDM}, \\
c\left( x\right) &\in \left[ 1,d_2 \right],~~ d_2 = const.>1,~ c(x) =1
\text{ for }x\in  \Omega _{FDM}, \\
a(x), c(x) &\in C^{2}\left( \mathbb{R}^{3}\right) . 
\end{split}
\end{equation}

\section{The domain decomposition algorithm}
\label{sec:hyb}

\begin{figure}  
 {\includegraphics[clip,scale=0.55, trim = 4.0cm 0.0cm 2.0cm 0.0cm, clip=true,]{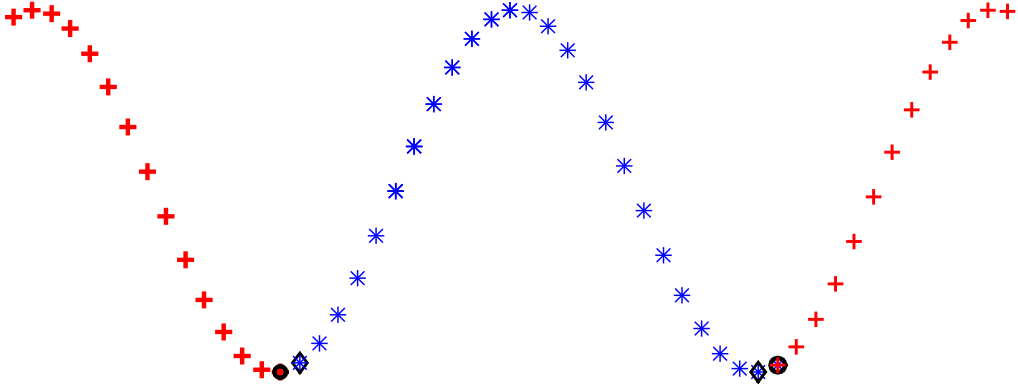}}
 \caption{Example of the solution in the domain decomposition between $\Omega_{FEM}$ and $\Omega_{FDM}$ in one
   dimension with $a=c=1$. The interior nodes of the  finite element
   grid  $\Omega_{FEM}$  are denoted by stars, while circles and diamonds denote nodes,
   which are shared between meshes in $\Omega_{FEM}$ and
   $\Omega_{FDM}$.  The circles are interior nodes  of the grid in
   $\Omega_{FDM}$, while the  diamonds are interior nodes  of the grid in
   $\Omega_{FEM}$. At each time iteration, solution obtained in
   $\Omega_{FDM}$ at $\omega_{\rm o}$ is copied to the corresponding boundary nodes in
   $\Omega_{FEM}$, while simultaneously the solution obtained in
   $\Omega_{FEM}$ at  $\omega_{\diamond}$ is copied to the corresponding boundary nodes in
   $\Omega_{FDM}$.}
 \label{fig:hyb1D}
\end{figure}

We now describe the domain decomposition between $\Omega_{FEM}$ and $\Omega_{FDM}$
domains.
 This communication is achieved by
mesh overlapping across a two-element thick layer around
$\Omega_{FEM}$ - see Figure ~\ref{fig:hyb1D}.
 First, using the Figure \ref{fig:hyb1D} we observe that the interior nodes of
the computational domain $\Omega$ belong to either of the following sets:
\begin{itemize}
\item[$\omega_{\rm o}$] Nodes '$\rm o$' - lie on the boundary of $\Omega_{FEM}$ and are interior to $\Omega_{FDM}$,
\item[$\omega_{\diamond}$] Nodes '$\diamond$' -  lie on the boundary of
$\Omega_{FDM}$ and are 
 interior to $\Omega_{FEM}$,
\item[$\omega_{*}$] Nodes '$*$' are interior to $\Omega_{FEM}$,
\item[$\omega_{+}$] Nodes '$+$' are interior to $\Omega_{FDM}$ 
\end{itemize}

Then the main loop in time for the explicit schemes which solves the  problem (\ref{modelhyb}) is as follows:

\textbf{Algorithm 1}

 At every time step $k$ we perform the
following operations:
\begin{enumerate}
\item On the structured part of the mesh $\Omega_{FDM}$ update the
 FDM  solution  at nodes $\omega_{+}$  and  $\omega_{\diamond}$.

\item On the unstructured part of the mesh $\Omega_{FEM}$ update the
  FEM solution at nodes $\omega_{*}$ and $\omega_{\rm o}$.

\item Copy FEM solution obtained at nodes $\omega_{\diamond}$ as a
  boundary condition for the FDM solution in $\Omega_{FDM}$.

\item Copy FDM solution obtained at nodes $\omega_{\rm o}$ as a
  boundary condition for the FEM solution in $\Omega_{FEM}$.

\end{enumerate}

\section{Energy estimate}

\label{sec:energyerror1}

 In this section we prove the uniqueness theorem, or energy estimate,
 for the function $u\in H^{2}\left( \Omega_{T}\right) $ of
  the problem (\ref{modelhyb}),
 using the technique of \cite{lad}.

 \textbf{Theorem}

 \emph{Assume that condition (\ref{coefic}) on the functions $c(x),
   a(x)$ hold. Let }$\Omega \subset \mathbb{R}^{3}$\emph{\ be a
   bounded domain with the piecewise smooth boundary }$\partial
 \Omega$. \emph{\ For any }$ t\in \left( 0,T\right) $\emph{\ let
 }$\Omega_{t}= \Omega\times \left( 0,t\right) .$ \emph{\ Suppose that
   there exists a solution }$u\in H^{2}(\Omega_{T})
 $\emph{\ of the problem (\ref{modelhyb}).  Then the function $u$ is
   unique and there exists a constant}
   $A=A(||c||_{\Omega},||a||_{\Omega}, t)$ \emph{\ such that the
   following energy estimate is true for all $c, a \geq 1$ in
   (\ref{modelhyb})}

\begin{equation}\label{estimate1}
\begin{split}
 \left \Vert \sqrt{\frac{1}{c^2}}\partial _{t} u(x,t)
 \right\Vert_{L_{2}\left( \Omega \right) }^{2} 
+  \left \Vert \sqrt{a} \nabla u (x,t) \right\Vert _{L_{2}\left( \Omega \right) }^{2}
  \leq A \left[ \left\Vert g \right\Vert _{L_{2}\left( \Omega _{t}\right)}^{2}+
\left\Vert \sqrt{\frac{1}{c^2} }f_{1}\right\Vert _{L_{2}\left( \Omega
\right) }^{2}+
\left\Vert \sqrt{a} \nabla f_{0}\right\Vert _{L_{2}\left( \Omega\right) }^{2}\right] .
\end{split}
\end{equation}

\textbf{Proof.}

First we multiply hyperbolic equation in (\ref{modelhyb}) by $2 \partial_t u$ and integrate
over $ \Omega\times \left( 0,t\right)$ to get
\begin{equation}\label{eq1_mod4}
\begin{split}
\int \limits_{0}^{t}\int\limits_{\Omega} 2~\frac{1}{c^2}
\partial_{tt} u ~\partial_t u~ dxd\tau 
- \int \limits_{0}^{t}\int\limits_{\Omega} 2 
\nabla \cdot (a \nabla u)~\partial_t u~ dxd\tau =
2\int\limits_{0}^{t}\int\limits_{\Omega} g~\partial_t u~dxd\tau.
\end{split}
\end{equation}

Integrating in time the first term of (\ref{eq1_mod4}) we get
\begin{equation}\label{eq4_time1}
\int\limits_{0}^{t}\int\limits_{\Omega}\partial _{t}( \frac{1}{c^2}
\partial_t u^{2})  dxd\tau 
=\int\limits_{\Omega}\left( \frac{1}{c^2}\partial_t u^{2} \right) \left( x,t\right) dx  
-\int\limits_{\Omega} \frac{1}{c^2} f_{1}^{2} \left( x,t\right) dx.
\end{equation}

Integrating by parts in space the second term of (\ref{eq1_mod4}),
using conditions (\ref{coefic}) giving $a=1$ on $\partial  \Omega$, and absorbing boundary conditions in
(\ref{modelhyb}) we get
\begin{equation}\label{grad1}
\begin{split}
&2\int\limits_{0}^{t}\int\limits_{\Omega} \nabla \cdot \left(a  \nabla  u \right) \partial_t u dx d\tau \\ 
&=
  2\int\limits_{0}^{t}\int\limits_{\partial \Omega}\left( \partial_t u
  \right) \partial_n u dS d\tau -
  2  \int\limits_{0}^{t}\int\limits_{\Omega}\left( a \nabla u \right) \left(
  \nabla \partial_t u \right) dxd\tau \\ &=
 -   \int\limits_{0}^{t}\int\limits_{\partial \Omega}\left( \partial_t u
  \right)^2 ~dS d\tau - \int\limits_{0}^{t}\int\limits_{\Omega} a \partial _{t}|\nabla u
  |^{2}dxd\tau.
\end{split}
\end{equation}

 Integrating last term of (\ref{grad1}) in time and using initial
 conditions of the equation (\ref{modelhyb}) we obtain
\begin{equation}\label{grad1_1}
\begin{split}
&-\int\limits_{0}^{t}\int\limits_{\Omega} a \partial _{t} |\nabla  u| ^{2}dxd\tau  
=-\int\limits_{\Omega} a | \nabla  u| ^{2}\left( x,t\right) dx
 + \int\limits_{\Omega} a |\nabla u |^{2}\left( x,0\right) dx   \\
&=-\int\limits_{\Omega}  a | \nabla  u |^{2}\left( x,t\right) dx
 + \int\limits_{\Omega}  a |\nabla  f_{0}| ^{2}\left( x\right) dx. 
\end{split}
\end{equation}

Next,
 collecting estimates (\ref{eq4_time1}),
(\ref{grad1}), (\ref{grad1_1}), using the fact that 
$ \int\limits_{0}^{t}\int\limits_{\partial \Omega}\left( \partial_t u \right)^2 ~dS d\tau \geq 0$
 and
substituting them in (\ref{eq1_mod4}) we have
\begin{equation}\label{eq_main}
\begin{split}
&\int\limits_{\Omega}\left( \frac{1}{c^2} \partial_t u^{2} \right) \left( x,t\right) dx  
+\int\limits_{\Omega} a | \nabla  u |^{2}\left( x,t\right) dx
 \\
&
\leq 2\int\limits_{0}^{t}\int\limits_{\Omega} |g|~|\partial_t u|~dxd\tau
+\int\limits_{\Omega} \frac{1}{c^2} f_{1}^{2} \left( x,t\right) dx
 + \int\limits_{\Omega} a |\nabla  f_{0}| ^{2}\left( x\right) dx.
\end{split}
\end{equation}

Finally, to estimate the first term in the right hand side of
(\ref{eq_main}) we use the arithmetic-geometric mean inequality $2ab
\leq a^2 + b^2$ to obtain
\begin{equation}\label{2ab}
2\int\limits_{0}^{t}\int\limits_{\Omega}| g |\cdot |\partial_t u|~dxd\tau \leq
\int\limits_{0}^{t}\int\limits_{\Omega} |g|^{2}dxd\tau +
\int\limits_{0}^{t}\int\limits_{\Omega}|\partial_t u|^{2}dxd\tau .  
\end{equation}

Noting that using (\ref{coefic})  we can write following estimate
\begin{equation*}
  \int\limits_{0}^{t}\int\limits_{\Omega} |\partial_t u|^{2} ~dxd\tau
\leq   \int\limits_{0}^{t}\int\limits_{\Omega} \left (\frac{1}{c^2}|\partial_t u|^{2} + a |\nabla u|^2 \right)
~dxd\tau,
\end{equation*}
and substituting  the above estimate into (\ref{2ab}) and then the resulting estimate into
(\ref{eq_main}) we have the following estimate
\begin{equation}\label{mod4_3}
\begin{split}
\int\limits_{\Omega} \left ( \frac{1}{c^2} \partial_t u^{2} + a\left| \nabla u \right|^{2} \right )(x,t) dx 
  &\leq \int\limits_{0}^{t}\int\limits_{\Omega} | g|^{2}dxd\tau 
 +   \int\limits_{0}^{t}\int\limits_{\Omega} \left ( \frac{1}{c^2}|\partial_t u|^{2} + a |\nabla u|^2 \right) ~dxd\tau \\
&+\int\limits_{\Omega}\left(  \frac{1}{c^2}  f_{1}^{2}+ a \left| \nabla f_{0}\right|^{2}\right)(x,t)dx.
\end{split}
\end{equation}

Let us denote 
\begin{equation}
F(t):=\int\limits_{\Omega} \left (  \frac{1}{c^2} \partial_t u^{2} + a \left| \nabla u\right|^{2} \right) \left( x,t\right) dx.  
\end{equation}
We rewrite estimate (\ref{mod4_3}) in the form
\begin{equation}\label{gronwall}
F(t) \leq A \int_0^t F(\tau) d\tau  + r(t),
\end{equation}
where $r(t):= \int\limits_{0}^{t}\int\limits_{\Omega} |g|^{2}dxd\tau 
 +\int\limits_{\Omega}\left(  \frac{1}{c^2} f_{1}^{2}+ a \left| \nabla f_{0}\right|
^{2}\right) \left( x,t\right) dx.$

Applying Gr\"onwall's inequality to (\ref{gronwall}) with a 
constant  $A=A(||c||_{\Omega},||a||_{\Omega}, t)$ we get desired energy estimate  (\ref{estimate1}) which also can be written in the form
\begin{equation}\label{mod4_5}
\begin{split}
\int\limits_{\Omega} \left ( \frac{1}{c^2} \partial_t u^{2} + a \left| \nabla u\right|^{2} \right)(x,t) dx
\leq
A \left( \int\limits_{0}^{t}\int\limits_{\Omega} |g|^{2}dxd\tau 
 +\int\limits_{\Omega}\left(  \frac{1}{c^2} f_{1}^{2}+ a \left| \nabla f_{0}\right|    
^{2} \right) \left( x,t\right) dx \right).
\end{split}
\end{equation}

$\square$

\section{Statement of the forward and inverse problems}

\label{sec:model}

In this section we state the forward and inverse problems. In section
\ref{subsec:ad_alg} we will show how these problems can be solved using the domain
decomposition algorithm of section \ref{sec:hyb}.

 Let the boundary $\partial
\Omega$ is such that $\partial \Omega =\partial _{1} \Omega \cup
\partial _{2} \Omega \cup \partial _{3} \Omega$ where $\partial _{1}
\Omega$ and $\partial _{2} \Omega$ are, respectively, front and back
sides of the domain $\Omega$, and $\partial _{3} \Omega$ is the union
of left, right, top and bottom sides of this domain.
Let at $S_T := \partial_1 \Omega \times (0,T)$ 
we will have time-dependent observations at the
backscattering side $\partial_1 \Omega$ of the domain $\Omega$.
We also define  $S_{1,1} := \partial_1
 \Omega \times (0,t_1]$, $S_{1,2} := \partial_1 \Omega \times
 (t_1,T)$,  $S_2 := \partial_2 \Omega \times (0, T)$ and  $S_3 :=
 \partial_3 \Omega \times (0, T)$.

We have used 2 model problems in our computations.

\textbf{Model Problem 1}

The first model problem  is the same as (\ref{modelhyb}) but when $c(x)=1 ~~\forall x \in \Omega$  and with non-homogeneous initial conditions:
\begin{equation}\label{model1}
\begin{split}
\frac{\partial^2 u}{\partial t^2}  -  \nabla \cdot ( a \nabla  u)   &= 0,~ \mbox{in}~~ \Omega_T, \\
  u(x,0) = f_0(x), ~~~u_t(x,0) &= f_1(x)~ \mbox{in}~~ \Omega,     \\
\partial _{n} u& = p\left( t\right) ~\mbox{on}~ S_{1,1},
\\
\partial _{n}  u& =-\partial _{t} u ~\mbox{on}~ S_{1,2},
\\
\partial _{n} u& =-\partial _{t} u~\mbox{on}~ S_2, \\
\partial _{n} u& =0~\mbox{on}~ S_3.\\
\end{split}
\end{equation}

\textbf{Model Problem 2}

Our model problem 2 uses homogeneous initial conditions and $c(x)=1 ~~\forall x \in \Omega$ in
(\ref{modelhyb})  and is defined
as
\begin{equation}\label{model2}
\begin{split}
\frac{\partial^2 u}{\partial t^2}  -  \nabla \cdot ( a\nabla  u)   &= 0,~ \mbox{in}~~ \Omega_T, \\
  u(x,0) = 0, ~~~u_t(x,0) &= 0~ \mbox{in}~~ \Omega,     \\
\partial _{n} u& = p\left( t\right) ~\mbox{on}~ S_{1,1},
\\
\partial _{n}  u& =-\partial _{t} u~\mbox{on}~ S_{1,2},
\\
\partial _{n} u& =-\partial _{t} u~\mbox{on}~ S_2, \\
\partial _{n} u& =0~\mbox{on}~ S_3.\\
\end{split}
\end{equation}

We assume that our
coefficient $a \left( x\right)$ of problems (\ref{model1}) and (\ref{model2})  is such
that
\begin{equation} \label{2.3}
\begin{split}
a \left( x\right) &\in \left[ 1,d \right],~~ d = const.>1,~ a(x) =1
\text{ for }x\in  \Omega _{FDM}, \\
a \left( x\right) &\in C^{2}\left( \mathbb{R}^{3}\right) . 
\end{split}
\end{equation}

We consider the following

\textbf{Inverse Problem 1 (IP1)} \emph{Suppose that the coefficient
}$a\left( x\right)$ \emph{\ in the problem (\ref{model1}) satisfies
  (\ref{2.3}).  Assume that the function }$ a\left( x\right)
$\emph{\ is unknown in the domain }$\Omega \diagdown
\Omega_{FDM}$\emph{. Determine the function }$ a\left( x\right)
$\emph{\  in (\ref{model1}) for }$x\in \Omega \diagdown \Omega_{FDM},$ \emph{\ assuming
  that the following function }$\tilde u\left( x,t\right) $\emph{\ is
  known}
\begin{equation}
  u\left( x,t\right) = \tilde u \left( x,t\right) ,\forall \left( x,t\right) \in
  S_T.  \label{2.5_1}
\end{equation}

The question of stability and uniqueness of \textbf{IP1} is
addressed in the recent work \cite{CristofolLiSoc}.

\textbf{Inverse Problem 2 (IP2)} \emph{Suppose that the coefficient
}$a \left( x\right)$ \emph{\ in the problem (\ref{model2}) satisfies
  (\ref{2.3}).  Assume that the function }$ a\left( x\right)
$\emph{\ is unknown in the domain }$\Omega \diagdown
\Omega_{FDM}$\emph{. Determine the function }$ a\left( x\right)
$\emph{\  in (\ref{model2}) for }$x\in \Omega \diagdown \Omega_{FDM},$ \emph{\ assuming
  that the following function }$\tilde u\left( x,t\right) $\emph{\ is
  known}
\begin{equation}
  u\left( x,t\right) = \tilde u \left( x,t\right) ,\forall \left( x,t\right) \in
  S_T.  \label{2.5}
\end{equation}

\section{Optimization method}

\label{sec:opt}

In this section we will reformulate our inverse problem \textbf{IP1}
as an optimization problem to be able to  reconstruct the unknown 
function $a(x)$ in (\ref{model1})  with best fit to time and space domain observations
$\tilde u$, measured at a finite number of observation points on
$\partial_1 \Omega$.  Solution of \textbf{IP2} follows from the
solution of \textbf{IP1} by taking $f_0 =f_1=0$.

Our goal is to minimize the Tikhonov functional
\begin{equation}
F(u, a) = \frac{1}{2} \int_{S_T}(u - \tilde{u})^2 z_{\delta }(t) dxdt +
\frac{1}{2} \gamma  \int_{\Omega}(a -  a_0)^2~~ dx,
\label{functional}
\end{equation}
where $\tilde{u}$ is the observed field, $u$ satisfies the equations
(\ref{model1}), and thus depends on $a$, $a_{0}$ is the initial guess
for $a$, and $\gamma$ is the regularization parameter.  Here,
$z_{\delta }(t)$ is a cut-off function, which is introduced to ensure
that the compatibility conditions at $\overline{\Omega}_{T}\cap
\left\{ t=T\right\} $ for the adjoint problem (\ref{adjoint1}). The function $z_{\delta }$ can be
chosen similarly with \cite{BCN}.

For our theoretical investigations we
introduce the following spaces of real valued  functions
\begin{equation}\label{spaces}
\begin{split}
H_u^1 &:= \{ w \in H^1(\Omega_T):  w( \cdot , 0) = 0 \}, \\
H_{\lambda}^1 &:= \{ w \in  H^1(\Omega_T):  w( \cdot , T) = 0\},\\
U^{1} &=H_{u}^{1}(\Omega_T)\times H_{\lambda }^{1}(\Omega_T)\times C\left( \overline{\Omega}\right),\\
U^{0} &=L_{2}\left(\Omega_{T}\right) \times L_{2}\left(\Omega_{T}\right) \times
L_{2}\left( \Omega \right). 
\end{split}
\end{equation}

To solve this minimization problem for model problem (\ref{model1})  we introduce the Lagrangian
\begin{equation}\label{lagrangian1}
\begin{split}
L(v) &= F(u, a) 
-  \int_{\Omega_T} \frac{\partial
 \lambda }{\partial t} \frac{\partial u}{\partial t}  ~dxdt  
+   \int_{\Omega_T}( a\nabla  u)( \nabla  \lambda)~dxdt  \\
&  - \int_{\Omega} \lambda(x,0) f_1(x) ~dx - \int_{S_{1,1}} \lambda p(t) ~d \sigma dt    + \int_{S_{1,2}} \lambda \partial_t u ~d\sigma dt
   + \int_{S_2} \lambda \partial_t u  ~d\sigma dt   , \\
\end{split}
\end{equation}
where $v=(u,\lambda, a) \in U^1$, and search for a stationary point
with respect to $v$ satisfying $ \forall \bar{v}= (\bar{u}, \bar{\lambda}, \bar{a}) \in U^1$
\begin{equation}
 L'(v; \bar{v}) = 0 ,  \label{scalar_lagr1}
\end{equation}
where $ L^\prime (v;\cdot )$ is the Jacobian of $L$ at $v$.

As usual, we assume that $\lambda \left( x,T\right) =\partial _{t}\lambda \left(
x,T\right) =0$ and  impose such conditions on the function
$\lambda $  that $ L\left( u,\lambda, a \right)
:=L\left( v\right) =F\left( u, a\right).$ 
We also  use
the facts that $\lambda (x ,T) = \frac{\partial \lambda}{\partial t}
(x,T) =0$  as
well as $a=1$ on $\partial \Omega$, together with initial conditions  of (\ref{model1}) and
boundary conditions $ \partial_n u = 0$  on $S_3$  and $  \partial_n \lambda = 0$ on
$\Omega_T \setminus S_T$.  The equation (\ref{scalar_lagr1}) expresses that
for all $\bar{v} \in U^1$,
\begin{equation}\label{forward1}
\begin{split}
0 = \frac{\partial L}{\partial \lambda}(u)(\bar{\lambda}) =
&- \int_{\Omega_T}  \frac{\partial \bar{\lambda}}{\partial t} \frac{\partial u}{\partial t}~ dxdt 
+  \int_{\Omega_T}  ( a \nabla u) (\nabla  \bar{\lambda}) ~ dxdt 
- \int_{\Omega} \bar{\lambda}(x,0) f_1(x) ~dx  \\
&- \int_{S_{1,1}} \bar{\lambda} p(t) ~d \sigma dt   + \int_{S_{1,2}} \bar{\lambda} \partial_t u ~d\sigma dt    + \int_{S_2} \bar{\lambda} \partial_t u  ~d\sigma dt,~~\forall \bar{\lambda} \in H_{\lambda}^1(\Omega_T),\\
\end{split}
\end{equation}
\begin{equation} \label{control1}
\begin{split}
0 = \frac{\partial L}{\partial u}(u)(\bar{u}) &=
\int_{S_T}(u - \tilde{u})~ \bar{u}~ z_{\delta}~ dxdt- \int_{\Omega} 
\frac{\partial{\lambda}}{\partial t}(0) \bar{u}(x,0) ~dx  \\
&-  \int_{\Omega_T}  \frac{\partial \lambda}{\partial t} \frac{\partial \bar{u}}{\partial t}~ dxdt
 + \int_{\Omega_T} ( a \nabla  \lambda) (\nabla  \bar{u})  ~ dxdt~~\forall \bar{u} \in H_{u}^1(\Omega_T).\\
\end{split}
\end{equation}
Finally, we obtain the equation which expresses that the gradient
with respect to  $a$  vanish:
\begin{equation} \label{grad1new} 
0 = \frac{\partial L}{\partial  a}(u)(\bar{a})
 =    \int_{\Omega_T} (\nabla  u) (\nabla  \lambda) \bar{a} ~dxdt 
+\gamma \int_{\Omega} (a - a_0) \bar{a}~dx,~ x \in \Omega.
\end{equation}

The equation (\ref{forward1}) is the weak formulation of the state equation
(\ref{model1}) and the equation (\ref{control1}) is the weak
formulation of the following adjoint problem
\begin{equation}
\begin{split} \label{adjoint1}
 \frac{\partial^2 \lambda}{\partial t^2} - 
  \nabla \cdot (a \nabla  \lambda)  &= -  (u - \tilde{u}) z_{\delta}, ~  x \in S_T,   \\
\lambda(\cdot, T)& =  \frac{\partial \lambda}{\partial t}(\cdot, T) = 0, \\
\partial _{n} \lambda& =0,~\mbox{on}~ \Omega_T\setminus S_T.
\end{split}
\end{equation}

We note that the Lagrangian (\ref{lagrangian1}) and the optimality
conditions (\ref{forward1}), (\ref{control1}) for the model problem 2
will be the same, as for the model problem 1, and only the difference will be
that  these expressions will note have terms with initial conditions.

\section{The domain decomposition FEM/FDM method}

\label{sec:domdec}

In this section we formulate finite element and finite difference
methods for the solution of model problem 2.  FEM for model problem 1
is the same only terms with non-zero initial conditions should be
induced in the discretization.

\subsection{Finite element discretization}
\label{sec:fem}

 We discretize $\Omega_{{FEM}_T} = \Omega_{FEM} \times (0,T)$ denoting by $K_h = \{K\}$ a partition of
 the domain $\Omega_{FEM}$ into tetrahedra $K$ ($h=h(x)$ being a mesh function,
 defined as $h |_K = h_K$, representing the local diameter of the elements),
 and we let $J_{\tau}$ be a partition of the time interval $(0,T)$ into time
 intervals $J=(t_{k-1},t_k]$ of uniform length $\tau = t_k - t_{k-1}$. We
 assume also a minimal angle condition on the $K_h$ \cite{Brenner}.


To formulate the finite element method,  we
 define the finite element spaces $C_h$, $W_h^u$ and $W_h^{\lambda}$.
First we introduce the finite element trial space $W_h^u$ for  $u$ defined by
\begin{equation}
W_h^u := \{ w \in H_u^1: w|_{K \times J} \in  P_1(K) \times P_1(J),  \forall K \in K_h,  \forall J \in J_{\tau} \}, \nonumber
\end{equation}
where $P_1(K)$ and $P_1(J)$ denote the set of linear functions on $K$
and $J$, respectively.
We also introduce the finite element test space  $W_h^{\lambda}$ defined by
\begin{equation}
W_h^{\lambda} := \{ w \in H_{\lambda}^1: w|_{K \times J} \in  P_1(K) \times P_1(J),  \forall K \in K_h,  \forall J \in J_{\tau} \}. \nonumber
\end{equation}

To approximate function
$a(x)$   we will use the space of piecewise constant functions $C_{h} \subset L_2(\Omega)$, 
\begin{equation}\label{p0}
C_{h}:=\{u\in L_{2}(\Omega ):u|_{K}\in P_{0}(K),\forall K\in  K_h\}, 
\end{equation}
where $P_{0}(K)$ is the piecewise constant function on $K$.

Next, we define $V_h = W_h^u \times W_h^{\lambda} \times C_h$.
 The finite element method now reads: Find $v_h \in V_h$, such
 that
\begin{equation}
L'(v_h)(\bar{v})=0 ~\forall
\bar{v} \in V_h .  \label{varlagr}
\end{equation}

The equation  (\ref{varlagr}) expresses that the finite element method for the forward problem (\ref{model2}) in $\Omega_{FEM}$ will be: Find $u_h \in W_h^u$, such that $\forall \bar{\lambda} \in W_h^\lambda$  and for known $a_h \in C_h$,
\begin{equation}\label{varforward}
\begin{split}
&- \int_{\Omega_{FEM_T}}  \frac{\partial \bar{\lambda}}{\partial t} \frac{\partial u_h}{\partial t}~ dxdt 
- \int_{\partial \Omega_{FEM}} \partial_n u_h \bar{\lambda} ~ dxdt
+  \int_{\Omega_{FEM_T}}  ( a_h \nabla u_h) (\nabla  \bar{\lambda}) ~ dxdt = 0,
\end{split}
\end{equation}
and the finite element method for the  adjoint problem (\ref{adjoint1}) in $\Omega_{FEM}$ reads: Find $\lambda_h \in W_h^\lambda$, such that $\forall \bar{u} \in W_h^u$ and for known $a_h \in C_h$,
\begin{equation} \label{varadjoint}
\begin{split}
&- \int_{\Omega_{FEM}} \frac{\partial{\lambda_h}}{\partial t}(0) \bar{u}(x,0) ~dx 
-  \int_{\Omega_{FEM_T}}  \frac{\partial \lambda_h}{\partial t} \frac{\partial \bar{u}}{\partial t}~ dxdt \\
&- \int_{\partial \Omega_{FEM}} \partial_n \lambda_h \bar{u} ~ dxdt
 + \int_{\Omega_{FEM_T}} ( a_h \nabla  \lambda_h) (\nabla  \bar{u})  ~ dxdt = 0.
\end{split}
\end{equation}
We note that usually $\dim V_{h}<\infty $ and $V_{h}\subset U^{1}$ as
a set and we consider $V_{h}$ as a discrete analogue of the space
$U^{1}.$ We introduce the same norm in $V_{h}$ as the one in $U^{0}$,
\begin{equation}\label{equiv}
\left\Vert \bullet \right\Vert _{V_{h}}:=\left\Vert \bullet \right\Vert
_{U^{0}},
\end{equation}
 where $U_0$ is defined in (\ref{spaces}).  From (\ref{equiv}) follows
 that in finite dimensional spaces all norms are equivalent. This
 allows us in numerical simulations compute coefficients in the space
 $C_h$. However, in the finite element discretization we write $a \in
 L_2(\Omega)$ to allow the function $a(x)$ be approximated in any other
 finite element space.

\subsection{Fully discrete scheme in $\Omega_{FEM}$}
\label{sec:discrete}

In this section we   present explicit schemes for computations of the solutions 
of forward and adjoint problems  in $\Omega_{FEM}$.
  After expanding functions $u_h(x,t)$ and $\lambda_h(x,t)$ in
terms of the standard continuous piecewise linear functions
$\{\varphi_i(x)\}_{i=1}^M$ in space and $\{\psi_k(t)\}_{k=1}^N$ in
time as 
\begin{equation*}
\begin{split}
u_h(x,t) &=\sum_{k=1}^N \sum_{i=1}^M u_{h_{i,k}}
\varphi_i(x)\psi_k(t), \\
\lambda_h(x,t) &=\sum_{k=1}^N \sum_{i=1}^M
\lambda_{h_{i,k}} \varphi_i(x)\psi_k(t), 
\end{split}
\end{equation*}
 where $u_{h_{i,k}}$ and
$\lambda_{h_{i,k}}$ denote unknown coefficients at the mesh point $x_i \in
K_h$ and time moment $t_k \in J_{\tau}$, substitute them into
(\ref{varforward}) and (\ref{varadjoint}), correspondingly, with
$\bar{\lambda}(x,t) =\bar{u}(x,t) =\sum_{l=1}^N \sum_{j=1}^M
\varphi_j(x)\psi_l(t)$. We note that we use finite element method only
inside $\Omega_{FEM}$, and thus we will have discrete solutions
$u_{h_{FDM}}:= u_{{h_{FDM}}_{i,k}}$ and $\lambda_{h_{FDM}} :=
\lambda_{{h_{FDM}}_{i,k}}$ obtained in $\Omega_{FDM}$ as the boundary
conditions at $\partial \Omega_{FEM}$, after exchange procedure.
We obtain the system of discrete
equations:
\begin{equation}\label{discrete1}
\begin{split}
&- \sum_{K \in \Omega_{FEM}} \sum_{k,l=1}^N \sum_{i,j=1}^M u_{h_{i,k}} \int_{K}
\varphi_i(x) \varphi_j(x) \int_{t_{k-1}}^{t_{k+1}} \partial_t \psi_k(t) \partial_t \psi_l(t)~dxdt \\
&- \sum_{\partial K \in \partial \Omega_{FEM}}  \sum_{k,l=1}^N \sum_{i,j=1}^M u_{h_{i,k}} \int_{\partial K}
\partial_n \varphi_i(x) \varphi_j(x) \int_{t_{k-1}}^{t_{k+1}}  \psi_k(t)  \psi_l(t)~dSdt \\
&+ \sum_{K \in \Omega_{FEM}} \sum_{k,l=1}^N \sum_{i,j=1}^M  u_{h_{i,k}} \int_{K}
a_h \nabla  \varphi_i(x) \nabla  \varphi_j(x) \int_{t_{k-1}}^{t_{k+1}}  \psi_k(t) \psi_l(t)~dxdt=0.
\end{split}
\end{equation}

For the case of adjoint problem (\ref{adjoint1}) we get the 
 system of  discrete equations:
\begin{equation}\label{discrete2}
\begin{split}
&- \sum_{K \in \Omega_{FEM}} \sum_{k,l=1}^N \sum_{i,j=1}^M \lambda_{h_{i,k}} \int_{K}
 \varphi_i(x) \varphi_j(x) \int_{t_{k-1}}^{t_{k+1}} \partial_t \psi_k(t) \partial_t \psi_l(t)~dxdt \\
&- \sum_{\partial K \in \partial \Omega_{FEM}}  \sum_{k,l=1}^N \sum_{i,j=1}^M 
 \lambda_{h_{i,k}} \int_{\partial K} \partial_n \varphi_i(x) \varphi_j(x) \int_{t_{k-1}}^{t_{k+1}} \psi_k(t)  \psi_l(t)~dSdt \\
&+ \sum_{K \in \Omega_{FEM}} \sum_{k,l=1}^N \sum_{i,j=1}^M \lambda_{h_{i,k}}  \int_{K} a_h
\nabla  \varphi_i(x) \nabla  \varphi_j(x) \int_{t_{k-1}}^{t_{k+1}}  \psi_k(t) \psi_l(t)~dxdt = 0.\\
\end{split}
\end{equation}  

Next, we compute explicitly time integrals in (\ref{discrete1}) and
(\ref{discrete2}) using the standard definition of piecewise-linear
functions in time, see \cite{BMaxwell} for details of this computation, and get the following systems of linear equations:
\begin{equation} \label{femod1}
\begin{split}
 M (\mathbf{u}^{k+1} - 2 \mathbf{u}^k  + \mathbf{u}^{k-1})  &=   - \tau^2 G (\frac{1}{6}\mathbf{u}^{k-1} 
+\frac{2}{3}\mathbf{u}^{k} + \frac{1}{6}\mathbf{u}^{k+1})  
 - \tau^2  K (\frac{1}{6}\mathbf{u}^{k-1} 
+\frac{2}{3}\mathbf{u}^{k} + \frac{1}{6}\mathbf{u}^{k+1}),   \\
M (\boldsymbol{\lambda}^{k+1} - 2 \boldsymbol{\lambda}^k + \boldsymbol{
\lambda}^{k-1}) &=  - \tau^2 G  (\frac{1}{6} \boldsymbol{\lambda}^{k-1} +\frac{2}{3} \boldsymbol{\lambda}^{k} + \frac{1}{6} \boldsymbol{\lambda}^{k+1})
 - \tau^2  K  (\frac{1}{6}  \boldsymbol{\lambda}^{k-1} +\frac{2}{3} \boldsymbol{\lambda}^{k} + \frac{1}{6} \boldsymbol{\lambda}^{k+1}), \\
\end{split}
\end{equation}
with initial  conditions :
\begin{eqnarray}
u(\cdot, 0)&= \frac{\partial u}{\partial t}(\cdot, 0) = 0, \\
\lambda (\cdot,T) &= \frac{\partial \lambda}{\partial t} (\cdot,T) =0.
\end{eqnarray}
  Here, $M$ is the block mass matrix in space, $K$ is the block stiffness
matrix, $G$ is the block matrix in space at  $\partial \Omega_{FEM}$,  $S^k$ is the load vector at
time level $t_k$, $\mathbf{u}^k$ and $ \boldsymbol{\lambda}^k$ denote the
nodal values of $u_h(\cdot,t_k)$ and $\lambda_h(\cdot,t_k)$, respectively, $\tau$
is the time step.

Now we define the mapping $F_K$ for the reference element $\hat{K}$
such that $F_K(\hat{K})=K$ and let $\hat{\varphi}$ be the piecewise
linear local basis function on the reference element $\hat{K}$ such
that $\varphi \circ F_K = \hat{\varphi}$.  Then the explicit formulas
for the entries in system (\ref{femod1}) at each element $K$ can be
given as:
\begin{equation}
\begin{split}
  M_{i,j}^{K} & =    ( ~\varphi_i \circ F_K, \varphi_j \circ F_K)_K, \\
  K_{i,j}^{K} & =   ( a_h \nabla  \varphi_i \circ F_K, \nabla  \varphi_j \circ F_K)_K,\\
 G_{i,j}^{\partial K} & = (\partial_n \varphi_i \circ F_K, \varphi_j \circ F_K)_{\partial K},\\
\end{split}
 \end{equation}
where $(\cdot,\cdot)_K$ denotes the $L_2(K)$ scalar product and $\partial K$ is the part of the boundary of element $K$ which lies at $\partial \Omega_{FEM}$.

To obtain an explicit scheme we approximate $M$ with the lumped mass
matrix $M^{L}$, see \cite{Cohen} for details about mass lumping
procedure.  We also   approximate terms corresponding to the mass  matrix in time, $\frac{1}{6}\mathbf{u}^{k-1} 
+\frac{2}{3}\mathbf{u}^{k} + \frac{1}{6}\mathbf{u}^{k+1}$  and $\frac{1}{6} \boldsymbol{\lambda}^{k-1} +\frac{2}{3} \boldsymbol{\lambda}^{k} + \frac{1}{6} \boldsymbol{\lambda}^{k+1}$, by $\mathbf{u}^{k}$ and $ \boldsymbol{\lambda}^{k}$, respectively, which fits to the procedure of  mass lumping in time.
Next, we multiply (\ref{femod1}) with $(M^{L})^{-1}$ and
get the following explicit method:
\begin{equation}   \label{fem}
\begin{split}
  \mathbf{u}^{k+1} = & (2 - \tau^2  (M^{L})^{-1} K)\mathbf{u}^k   - \tau^2  (M^{L})^{-1} G \mathbf{u}^k    -\mathbf{u}^{k-1},  \\
  \boldsymbol{\lambda}^{k-1} = &
  - \tau^2  (M^{L})^{-1} G \boldsymbol{\lambda}^k   + (2  - \tau^2  (M^{L})^{-1} K) \boldsymbol{\lambda}^k  -\boldsymbol{\lambda}^{k+1}.  \\
\end{split}
\end{equation} 

Finally, for reconstructing $a(x)$ in $\Omega_{FEM}$ we can use a gradient-based method
with an appropriate initial guess values of $a_0$.  The
discrete version in space of the gradient with respect to coefficient
$a$  in (\ref{grad1new}) take the form:
\begin{equation} \label{gradient1}
g_h(x) =   \int_{0}^T \nabla u_h \nabla  \lambda_h dt + \gamma (a_h - a_0).
\end{equation}

Here, $\lambda_h$  and $u_h$ are computed values of the adjoint
and forward problems  using explicit schemes
(\ref{fem}), and $a_h$   are approximated values of the computed
coefficient.

  \subsection{Finite difference formulation}
\label{sec:fdm}

 We recall now that from  conditions (\ref{2.3})  it follows that in $\Omega_{FDM}$
 the function $a(x)= 1$. This means that in
 $\Omega_{FDM}$ for the model problem (\ref{model1}) the forward problem will be
\begin{equation}\label{FDMmodel2}
\begin{split}
\frac{\partial^2 u}{\partial t^2}  -  \Delta  u   &= 0,~ \mbox{in}~~ \Omega_T, \\
  u(x,0) = f_0(x), ~~~u_t(x,0) &= f_1(x)~ \mbox{in}~~ \Omega_{FDM},     \\
\partial _{n} u& = p\left( t\right) ~\mbox{on}~ S_{1,1},
\\
\partial _{n}  u& =-\partial _{t} u~\mbox{on}~ S_{1,2},
\\
\partial _{n} u& =-\partial _{t} u~\mbox{on}~ S_2, \\
\partial _{n} u& =0~\mbox{on}~ S_3,\\
\partial _{n} u& = \partial _{n} u_{FEM}~\mbox{on}~ \partial \Omega_{FEM}.\\
\end{split}
\end{equation}
Then the corresponding  adjoint problem in $\Omega_{FDM}$ will be:
\begin{equation}
\begin{split} \label{FDMadjoint}
 \frac{\partial^2 \lambda}{\partial t^2} - 
  \Delta   \lambda  &= -  (u - \tilde{u}) z_{\delta}, ~  x \in S_T,   \\
\lambda(\cdot, T)& =  \frac{\partial \lambda}{\partial t}(\cdot, T) = 0, \\
\partial _{n} \lambda& =0,~\mbox{on}~ \Omega_T\setminus S_T, \\
\partial _{n} \lambda & = \partial _{n} \lambda_{{FEM}}~\mbox{on}~ \partial \Omega_{FEM}.\\
\end{split}
\end{equation}

Using standard finite difference discretization of the
equation~(\ref{FDMmodel2}) in $\Omega_{FDM}$ we obtain the following explicit
scheme for the solution of forward problem:
\begin{equation}  \label{fdmschemeforward}  
  u_{l,j,m}^{k+1} =  \tau^2 \Delta u_{l,j,m}^k + 2 u_{l,j,m}^k - u_{l,j,m}^{k-1},
\end{equation}
and the following explicit scheme for the adjoint problem which we solve backward in time:
\begin{equation}  \label{fdmschemeadjoint}  
  \lambda_{l,j,m}^{k-1} = -\tau^2 (u - \tilde{u})_{l,j,m}^k z_{\delta} 
 + \tau^2 \Delta \lambda_{l,j,m}^k + 2 \lambda_{l,j,m}^k - \lambda_{l,j,m}^{k+1},
\end{equation}
with  corresponding boundary conditions for every
problem.  In equations above, $u_{l,j,m}^{k}$ is the solution on the
time iteration $k$ at the discrete point $(l,j,m)$, 
 $(u - \tilde{u})_{{l,j,m}}^k$ is
the discrete analog of the difference $u - \tilde{u}$ at the
observations points at $S_T$, $\tau$ is the time step, and $ \Delta
u_{l,j,m}^k$ is the discrete Laplacian. In three dimensions, to
approximate $ \Delta u_{l,j,m}^k$ we get the standard seven-point
stencil:
\begin{eqnarray}
  \Delta u_{l,j,m}^k &=& \frac {u_{l+1,j,m}^k - 2 u_{l,j,m}^k 
    + u_{l-1,j,m}^k}{ \delta x_1^2 } + \frac {u_{l,j+1,m}^k - 2 u_{l,j,m}^k + 
u_{l,j-1,m}^k}{ \delta x_2^2 } + \nonumber \\ 
&&  
\frac {u_{l,j,m+1}^k - 2 u_{l,j,m}^k + u_{l,j,m-1}^k}{ \delta x_3^2 } ,
\end{eqnarray} 
where $\delta x_1$, $\delta x_2$, and $\delta x_3$ are the steps of the discrete finite 
difference meshes in the directions $x_1,x_2,x_3$, respectively.

\subsection{Absorbing boundary conditions}
\label{sec:absbc}

In our domain decomposition method we use first order absorbing
boundary conditions \cite{EM} which are exact for the case of our
computational tests of section \ref{sec:Numer-Simul}. We note that
these boundary conditions are implemented in efficient way in the
software package WavES \cite{waves} in the domain decomposition
method, and this is the main point of application of this method in
numerical simulations of section \ref{sec:Numer-Simul}.

To discretize first order absorbing boundary conditions \cite{EM}
\begin{equation}\label{absbc}
\frac{\partial u(x,t)}{\partial n} = -\frac{\partial u(x,t)}{\partial t} 
\end{equation}
 at the outer boundary of $\Omega_{FDM}$
we use forward finite difference approximation in the middle
point. This allows to obtain  
numerical approximation of higher order than   ordinary (backward or forward)
finite difference approximation.  If we  discretize the left boundary of $\Omega_{FDM}$ then 
we have  the condition (\ref{FDMmodel2})  in the form
\begin{equation*}
\frac{\partial u(x,t)}{\partial x_1} = \frac{\partial u(x,t)}{\partial t}.
\end{equation*}
Then the forward finite difference approximation in the middle
point of the above equation will be resulted in the following discretization
\begin{equation} \label{eq:abs1}
  \frac{ u_{l,j,m}^{k+1} - u_{l,j,m}^k}{\tau} +
  \frac{ u_{l+1,j,m}^{k+1} - u_{l+1,j,m}^k}{\tau} - \\
  \frac{ u_{l+1,j,m}^{k} - u_{l,j,m}^k}{ \delta x_1} - \frac{
    u_{l+1,j,m}^{k+1} - u_{l,j,m}^{k+1}}{ \delta x_1} = 0,
\end{equation}
which can be transformed to
\begin{equation} \label{eq:abs2}
  u_{l,j,m}^{k+1} = u_{l+1,j,m}^k + u_{l,j,m}^k \frac{ \delta x_1 - \tau}{ \delta x_1 + \tau}  -  \\
  u_{l+1,j,m}^{k+1}\frac{ \delta x_1- \tau}{ \delta x_1 + \tau},
\end{equation}
where  $ \delta x_1$  is the mesh size  in $x_1$ direction.
For other boundaries of the  $\Omega_{FDM}$ absorbing  boundary conditions
 can be written  similarly.

\subsection{The domain decomposition algorithm to solve forward and adjoint problems}
\label{sec:hybalg}

In this section we will present domain decomposition algorithm for the
solution of state and adjoint equations.
  We  note that because of
using explicit domain decomposition FEM/FDM method we need to choose
time step $\tau$ such that the whole scheme remains stable. We use the
stability analysis on the structured meshes \cite{Cohen} and choose the largest
time step in our computations accordingly to the CFL stability
condition  
\begin{equation} \label{CFL0}
\tau \leq \frac{1}{  \sqrt{a} \sqrt{\frac{1}{\delta x_1^2} + \frac{1}{\delta x_2^2} + \frac{1}{\delta x_3^2}}}. 
\end{equation}
Usually, we have the same mesh size $ \delta x_1=  \delta x_2= \delta x_3=h$  in $(x_1,x_2,x_3)$ directions, and the condition (\ref{CFL0}) can be
rewritten in three dimensions as
\begin{equation} \label{CFL}
\tau \leq h \sqrt{\frac{1}{ 3 a}}. 
\end{equation}

\textbf{Algorithm 2}

 At every time step $k$ we perform the
following operations:
\begin{enumerate}
\item On the structured part of the mesh $\Omega_{FDM}$ compute
  $u^{k+1}$, $\lambda^{k-1}$ from (\ref{fdmschemeforward}),
  (\ref{fdmschemeadjoint}), correspondingly, with $u^k, u^{k-1}$ and
  $\lambda^k, \lambda^{k+1}$ known.
\item On the unstructured part of the mesh $\Omega_{FEM}$  compute
$u^{k+1}, \lambda^{k-1}$ by using the explicit finite element schemes
(\ref{femod1}), correspondingly,  with $u^k, u^{k-1}$ and
  $\lambda^k, \lambda^{k+1}$ known.
\item Use the values of the function $u^{k+1}, \lambda^{k-1}$ at nodes
  $\omega_{\diamond}$, which are computed using the finite element
  schemes (\ref{femod1}), as a boundary condition for the finite
  difference method in $\Omega_{FDM}$.
\item Use the values of the functions $u^{k+1}, \lambda^{k-1}$ at
  nodes $\omega_{\rm o}$, which are computed using the finite difference
  schemes (\ref{fdmschemeforward}), (\ref{fdmschemeadjoint}),
  correspondingly, as a boundary condition for the finite element
  method in $\Omega_{FEM}$.
\item Apply swap of the solutions for the computed functions $u^{k+1},
  \lambda^{k-1}$ to be able perform algorithm on a new time level $k$.
\end{enumerate}

\section{The  algorithm for the solution of an inverse problem}
\label{subsec:ad_alg}

We use conjugate gradient method for iterative update of
approximations $a_{h}^{m}$ of the function $a_{h}$, where $m$ is the
number of iteration in our optimization procedure. We denote
\begin{equation}
\begin{split}
{g}^m(x) = \int_{0}^T \nabla  u_h^m \nabla  \lambda_h^m  dt + \gamma (a_h^m - a_0), \\
\end{split}
\end{equation}
where functions $u_{h}\left( x,t,a_{h}^{m}\right) ,\lambda _{h}\left(x,t,a_{h}^{m}\right) $\ are computed by solving
the state and adjoint problems
with $a:=a_{h}^{m}$.\\

\textbf{Algorithm  3}

\begin{itemize}
\item[Step 0.]  Choose the mesh $K_{h}$ in $\Omega$ and time partition $J_{\tau}$ of the time interval $\left( 0,T\right) .$
Start with the initial approximation $a_{h}^{0}= a_0$ and compute the
sequences of $a_{h}^{m}$ via the following steps:

\item[Step 1.]  Compute solutions $u_{h}\left( x,t,a_{h}^{m}\right) $
  and $\lambda _{h}\left( x,t,a_{h}^{m}\right) $ of state (\ref{model2}) and adjoint (\ref{adjoint1}) problems on $K_{h}$ and
  $J_{\tau}$ using domain decomposition algorithm of section \ref{sec:hybalg}.

\item[Step 2.]  Update the coefficient $a_h:=a_{h}^{m+1}$ 
  on $K_{h}$ and $J_{\tau}$ using the conjugate gradient method
\begin{equation*}
\begin{split}
a_h^{m+1} &=  a_h^{m}  + \alpha d^m(x),\\
\end{split}
\end{equation*}
where $\alpha$, is the step-size in the gradient update \cite{Peron} and
\begin{equation*}
\begin{split}
 d^m(x)&=  -g^m(x)  + \beta^m  d^{m-1}(x),
\end{split}
\end{equation*}
with
\begin{equation*}
\begin{split}
 \beta^m &= \frac{|| g^m(x)||^2}{|| g^{m-1}(x)||^2},
\end{split}
\end{equation*}
where $d^0(x)= -g^0(x)$.
\item[Step 3.]  Stop computing $a_{h}^{m}$ and obtain the function
  $a_h$ if either $||g^{m}||_{L_{2}( \Omega)}\leq \theta$ or norms
  $||g^{m}||_{L_{2}(\Omega)}$ are stabilized. Here, $\theta$ is the
  tolerance in updates $m$ of gradient method.  Otherwise set $m:=m+1$
  and go to step 1.

\end{itemize}

\section{Numerical Studies}
\label{sec:Numer-Simul}

\begin{figure}[tbp]
 \begin{center}
 \begin{tabular}{cc}
 {\includegraphics[scale=0.3, clip = true, trim = 0.0cm 0.0cm 0.0cm 0.0cm]{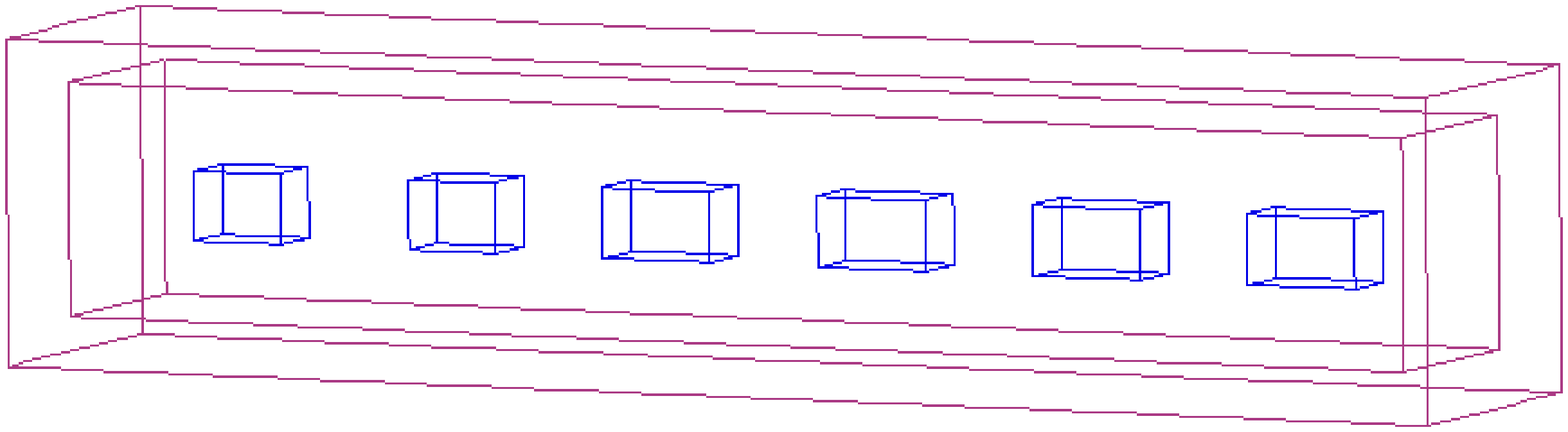}} &
 {\includegraphics[scale=0.3, clip = true, trim = 0.0cm 0.0cm 0.0cm 0.0cm]{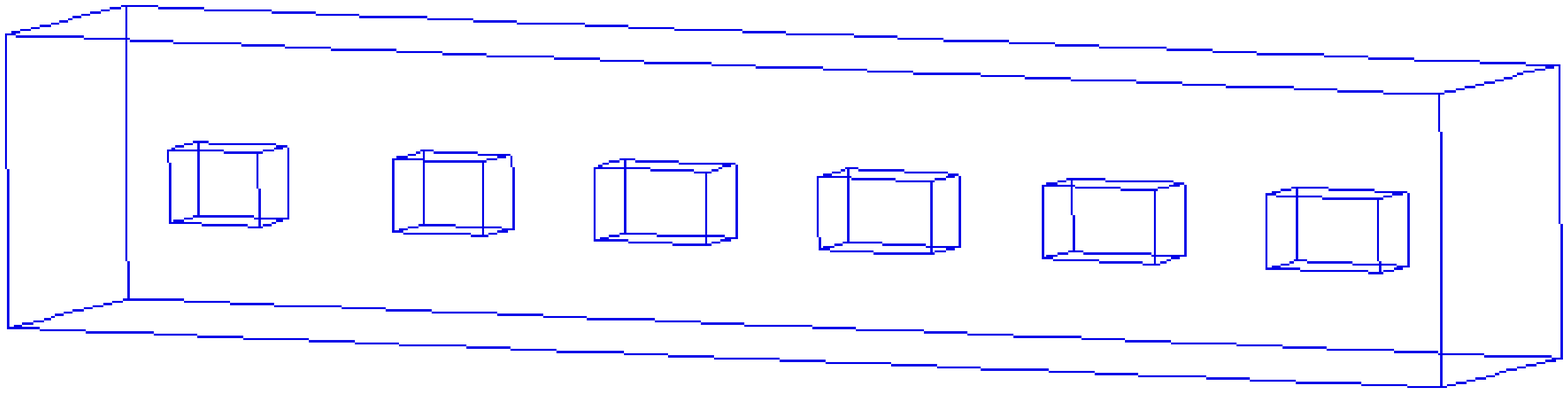}}
 \\
 a) $\Omega =  \Omega_{FEM} \cup \Omega_{FDM}$ &  b)   $\Omega_{FEM}$   \\
 \end{tabular}
 \end{center}
 \caption{{\protect\small \emph{ a) The hybrid domain  $\Omega= \Omega_{FEM} \cup \Omega_{FDM}$ b)  Finite element  domain $\Omega_{FEM}$.}}}
 \label{fig:fig1}
 \end{figure}

In this section we present numerical simulations of the reconstruction
of unknown function $a(x)$ inside a domain $\Omega_{FEM}$ using the
algorithm of section \ref{subsec:ad_alg}.  Accordingly to the
condition (\ref{2.3}) this function is known inside $\Omega_{FDM}$ and is set
to be $a(x) =1$.  The goal of our numerical tests is to reconstruct
scatterers of waveguide of figure \ref{fig:fig1} with $c=4.0$ inside
every small scatterer of figure \ref{fig:fig1}.

 The computational
 geometry $\Omega$ is split
 into two geometries, $\Omega_{FEM}$ and $\Omega_{FDM}$ such that
 $\Omega = \Omega_{FEM} \cup \Omega_{FDM}$, see Figure \ref{fig:fig1}. Next, we introduce dimensionless spatial variables 
 $x^{\prime}= x/\left(1m\right)$ and obtain that the domain
 $\Omega_{FEM}$ is transformed into  dimensionless
 computational domain
 \begin{equation*}
 \Omega_{FEM} = \left\{ x= (x_1,x_2,x_3) \in (
 -3.2,3.2) \times (-0.6,0.6) \times (-0.6,0.6) \right\} .
 \end{equation*}
   The dimensionless size of our computational domain
 $\Omega$ for the forward problem is
 \begin{equation*}
 \Omega = \left\{ x= (x_1,x_2,x_3) \in (
 -3.4,3.4) \times (-0.8,0.8) \times (-0.8,0.8) \right\} .
 \end{equation*}
 The space mesh in $\Omega_{FEM}$ and in $\Omega_{FDM}$ consists of
 tetrahedral and cubes, respectively. 
  We choose the mesh size $h=0.1$
 in our geometries in the hybrid FEM/FDM method, as well as in the
 overlapping regions between FEM and FDM domains.

In all our computations we use
  single  plane wave $p(t)$ initialized at $\partial_1 \Omega$
  in time $T=[0,3.0]$ such that
 \begin{equation}\label{f}
 \begin{split}
 p\left( t\right) =\left\{ 
 \begin{array}{ll}
 \sin \left( \omega t \right) ,\qquad &\text{ if }t\in \left( 0,\frac{2\pi }{\omega }
 \right) , \\ 
 0,&\text{ if } t>\frac{2\pi }{\omega }.
 \end{array}
 \right. 
 \end{split}
 \end{equation}

For generation of simulated backscattered data we define exact
function ${a(x)}=4$ inside small scatterers, see Figure
\ref{fig:fig1}, and $a(x)=1$ at all other points of the computational
domain $\Omega_{FEM}$.  Then we solve the forward problem on a refined mesh
which is not the same as used in our computations of inverse
problem. In a such way we avoid the problem with variational crimes.
The time step in all our computations is set to be  $\tau=0.006$ which
 satisfies the CFL condition \cite{CFL67}.
Isosurfaces of the simulated solution for the problem (\ref{model1})
with exact function $a(x)$ and $\omega=60$ in (\ref{f}) are presented
in figure \ref{fig:forward}. Using this figure we observe non-zero behavior of this
solution with initialized initial condition (\ref{initcond}).

In all our numerical simulations we have considered the additive noise
$\sigma$ introduced to the simulated boundary data $\tilde{u}$ in
(\ref{2.5}).
We have performed following reconstruction tests with the same values
of parameters in the reconstruction algorithm 3:

\begin{itemize}

\item  Test 1. Reconstruction of the function $a(x)$ in Model Problem 1 for $\omega \in [20,60]$ 
in (\ref{f}) and additive noise levels $\sigma= 3\%$  and  $\sigma= 10\%$.

\item  Test 2. Reconstruction of the function $a(x)$ in Model Problem 2 for $\omega \in [20,60]$ in (\ref{f}) and additive noise levels $\sigma= 3\%$  and  $\sigma= 10\%$.

\end{itemize}

In all our tests 
   we start the  algorithm 3 with guess values of the
   parameter $a(x)=1.0$ at all points in
   $\Omega$. We refer to  \cite{B}-\cite{BMaxwell}, \cite{HybEf, BJ}  for a similar choice of initial
   guess  which corresponds to starting of the algorithm 3 
 from the homogeneous
   domain.  The minimal and maximal values of the functions
   $a(x)$ in our computations belongs to the
   following set of admissible parameters
 \begin{equation}\label{admpar}
 \begin{split}
  M_{a} \in \{a\in C(\overline{\Omega })|1\leq a(x)\leq 5\}.\\
 \end{split}
 \end{equation}

We regularize 
the solution of the inverse problem
 by  computations with single value of the regularization
parameter $\gamma =0.01$ in (\ref{functional}). 
 Our computational experience have shown that such choice of $\gamma$
 is optimal one in our case.  Testing of different techniques, see , for example, \cite{Engl}, of the
 computing of regularization parameter is the topic of our ongoing
 research.
 The tolerance
 $\theta$ in our algorithm (section \ref{subsec:ad_alg}) is set to
 $\theta=10^{-6}$.

We use a post-processing procedure to get final images with our
reconstructions.  This procedure is as follows: assume, that functions
$a^m(x)$ are our reconstructions obtained by algorithm 3 of section
\ref{subsec:ad_alg} where $m$ is the number of iteration when we have
stopped to compute $a(x)$. Then to get post-processed images, we set
 \begin{equation}\label{postproc}
 \widetilde{a}^m(x)=\left\{ 
 \begin{array}{ll}
 a^m(x), & \text{ if }a^m(x)> 0.6 \max\limits_{\Omega_{FEM} }a^m(x), \\ 
 1, & \text{ otherwise. }
 \end{array}
 \right. 
 \end{equation}

Results of reconstruction for both tests are presented in tables 1,2.
Here, computational errors in procents are computed for $
\max\limits_{\Omega_{FEM} } a_{\overline{N}}$, where $\overline{N}
:=m$, and are compared with exact ones $a(x)=4$.

\begin{table}[h] 
{\footnotesize Table 1. \emph{Results of reconstruction
    of $a(x)$ for $\sigma =3\%$ together with computational errors in
    procents. Here, $\overline{N}$ is the final iteration number $m$
    in the conjugate gradient method of section \ref{subsec:ad_alg}.}}
\par
\vspace{2mm}
\centerline{
\begin{tabular}{|c|c|}
 \hline
   $\sigma= 3\%$ &  $\sigma = 3\%$ 
 \\
 \hline
\begin{tabular}{c|c|c|c} \hline
Test 1 & $ \max\limits_{\Omega_{FEM} } a_{\overline{N}}$ & error, $\%$  & $\overline{N}$  \\ \hline
$\omega=20$ & $5.0$ & 25 & $7  $   \\
$\omega=30$ & $5$ &  25   &$9$   \\
$\omega=40$ & $4.73$ & 18.3  &$10$   \\
$\omega=50$  & $5.0$ & 25 &$ 11$   \\
$\omega=60$  & $5.0$ & 25   &$12$   \\
\end{tabular}
 & 
\begin{tabular}{c|c|c|c} \hline
Test 2 & $ \max\limits_{\Omega_{FEM} } a_{\overline{N}}$ & error, $\%$  & $\overline{N}$ \\ \hline
$\omega=20$ & $5.0 $ & 25  & $ 8 $   \\
$\omega=30$ & $4.86$ & 21.5  & $9$   \\
$\omega=40$ & $4.19$ &  4.75   & $10$   \\
$\omega=50$  & $5.0$ &  25   & $11 $   \\
$\omega=60$  & $5.0$ & 25  & $12$   \\
\end{tabular} 
\\
\hline
\end{tabular}}
\end{table}

\begin{table}[h] 
{\footnotesize Table 2. \emph{Results of
    reconstruction of $a(x)$ for $\sigma =10\%$  together with computational errors in procents. Here, $\overline{N}$ is the final iteration number
    $m$ in the conjugate gradient method of section
    \ref{subsec:ad_alg}.}}  \par
\vspace{2mm}
\centerline{
\begin{tabular}{|c|c|}
 \hline
   $\sigma= 10\%$ &  $\sigma = 10\%$ 
 \\
 \hline
\begin{tabular}{c|c|c|c} 
Test 1 & $ \max\limits_{\Omega_{FEM} } a_{\overline{N}}$ & error, $\%$ & $\overline{N}$  \\
 \hline
$\omega=20$ & $3.29$ & 17.75 &  $ 8 $   \\
$\omega=30$ & $4.94$ & 23.5 & $10$   \\
$\omega=40$ & $4.35$ & 8.75  &  $11$   \\
$\omega=50$  & $4.4$ & 10 & $12 $   \\
$\omega=60$  & $4.44$ & 11 & $13$ \\  
\end{tabular}
 & 
\begin{tabular}{c|c|c|c}
Test 2 & $ \max\limits_{\Omega_{FEM} } a_{\overline{N}}$ & error, $\%$ & $\overline{N}$ \\ 
\hline
$\omega=20$ & $3.28 $ &  18  & $8$   \\
$\omega=30$ & $4.57$  & 14.25  & $10$   \\
$\omega=40$ & $3.96$  &  1  & $11$   \\
$\omega=50$  & $4.18$ &  4.5 & $ 12$   \\
$\omega=60$  & $4.46$ &  11.5 & $13 $   \\
\end{tabular} 
\\
\hline
\end{tabular}}
\end{table}

 \subsection{Test 1}\label{sec:casei}

 \begin{figure}[tbp]
 \begin{center}
 \begin{tabular}{cc}
 {\includegraphics[scale=0.3, angle=-90, trim = 1.0cm 1.0cm 1.0cm 1.0cm, clip=true,]{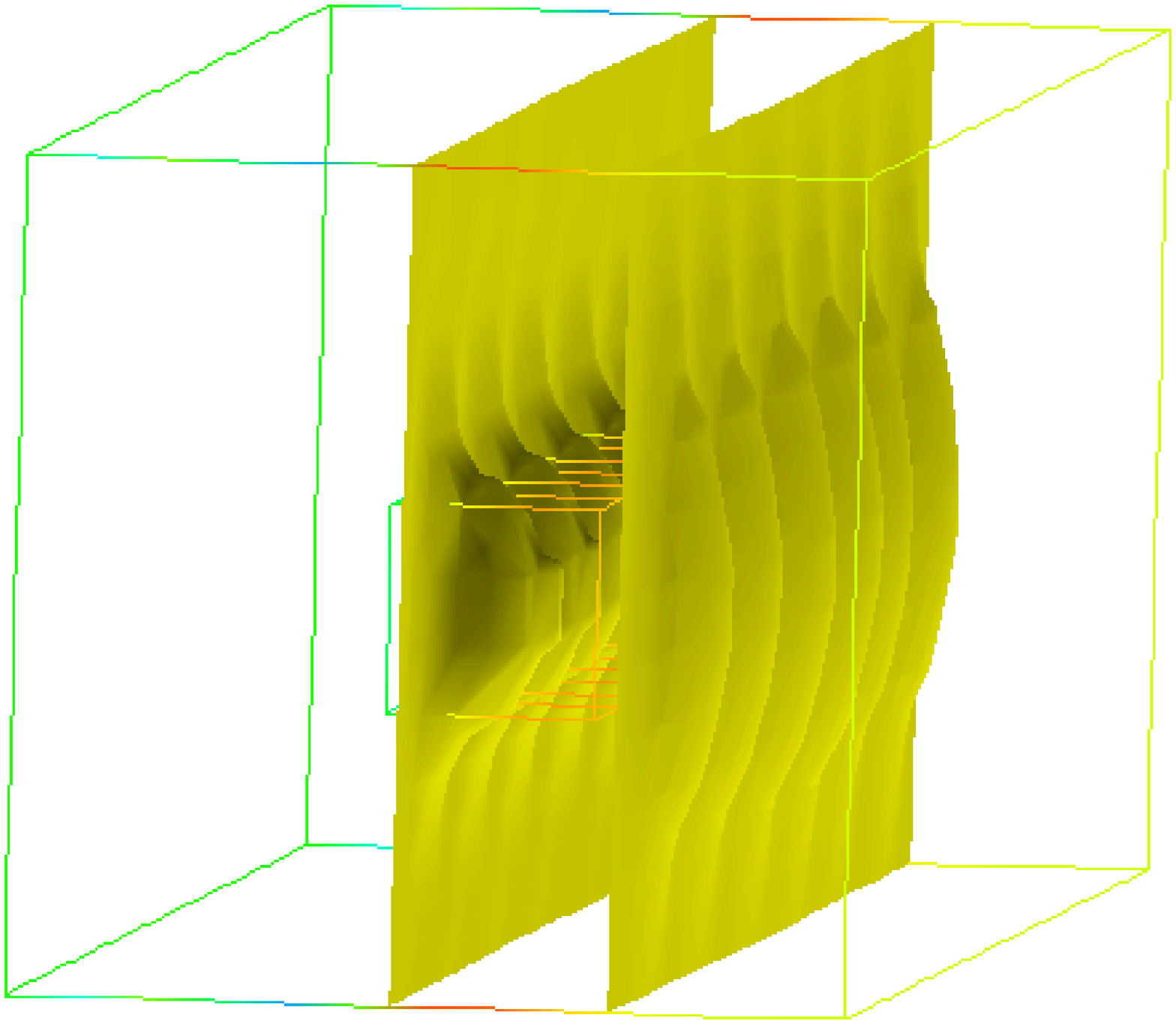}} &
 {\includegraphics[scale=0.35, angle=-90, trim = 1.0cm 1.0cm 1.0cm 1.0cm, clip=true,]{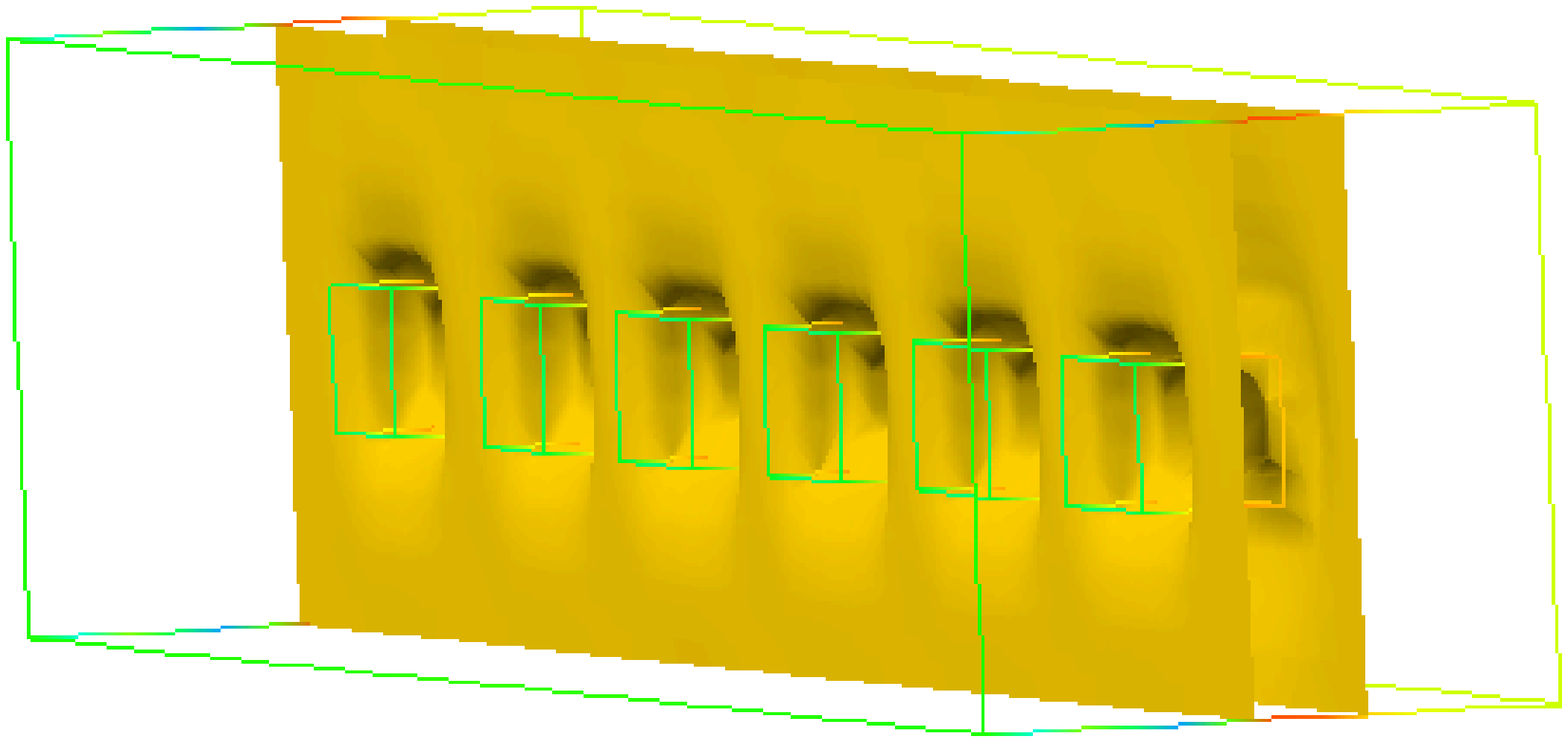}} \\
 a) $t=1.5$ & b) $t=1.5$ \\
 {\includegraphics[scale=0.3, angle=-90, trim = 1.0cm 1.0cm 1.0cm 1.0cm, clip=true,]{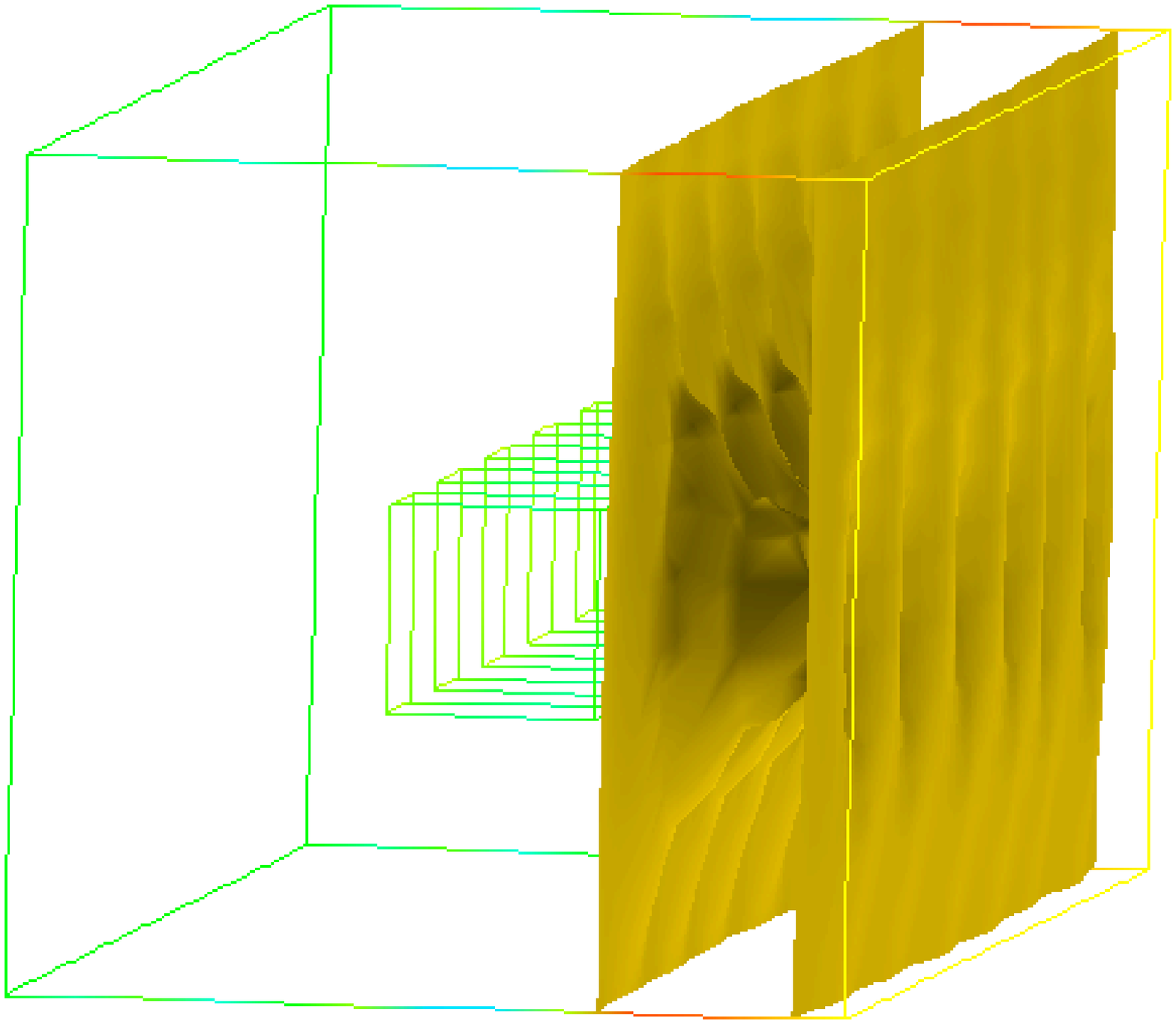}} &
 {\includegraphics[scale=0.35, angle=-90, trim = 1.0cm 1.0cm 1.0cm 1.0cm, clip=true,]{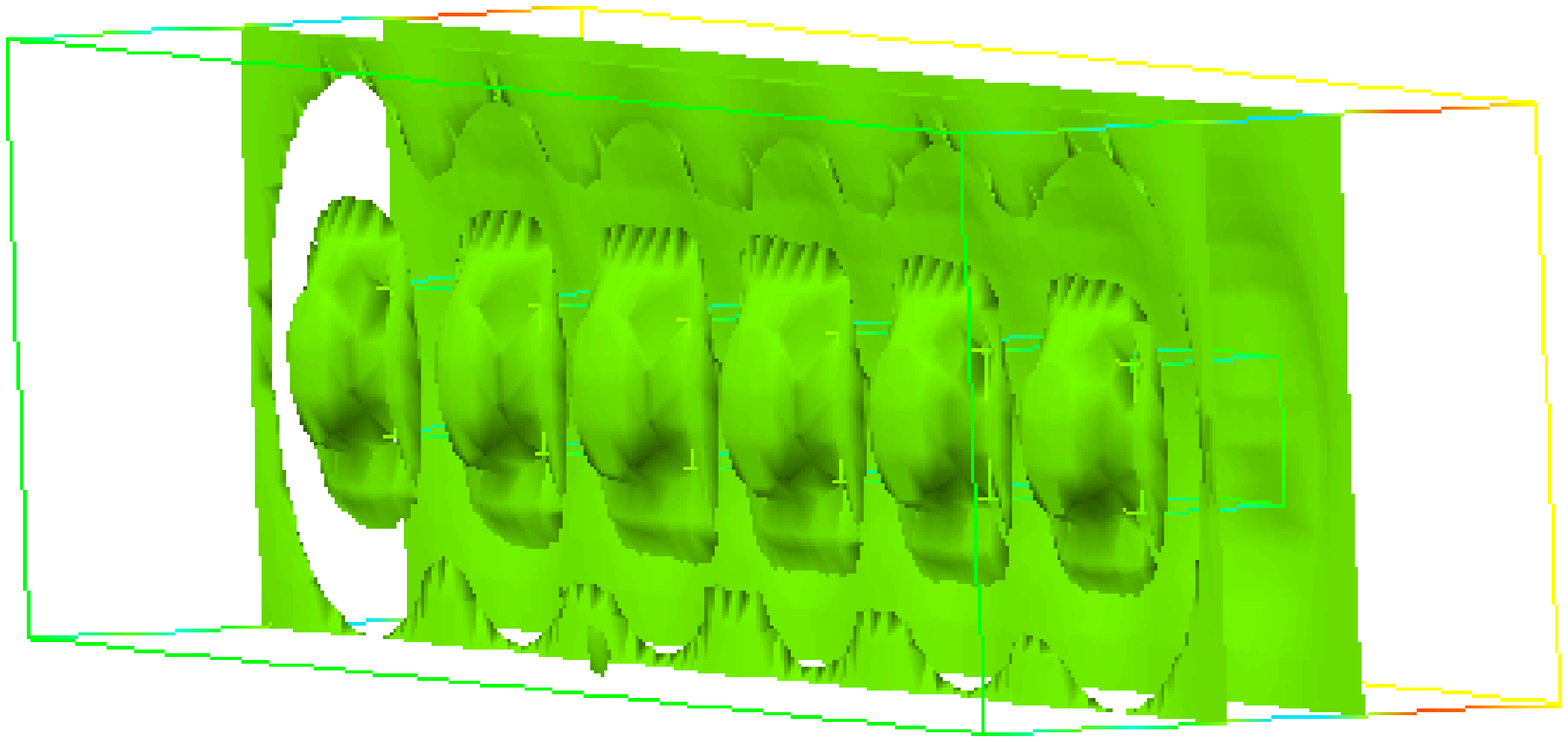}} \\
 c) $t=1.8$ & d) $t=1.8$  \\
 {\includegraphics[scale=0.3, angle=-90, trim = 1.0cm 1.0cm 1.0cm 1.0cm, clip=true,]{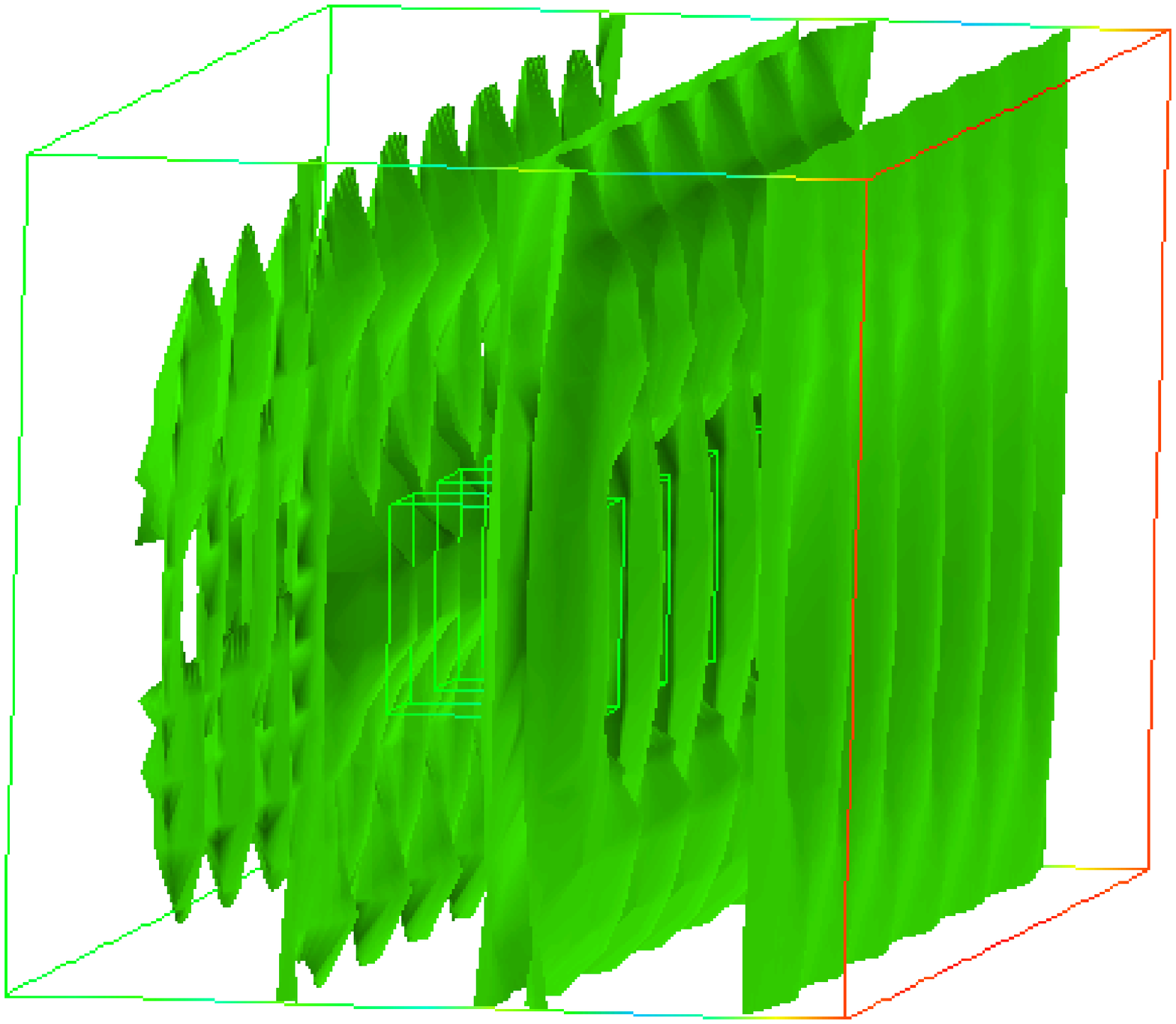}} &
 {\includegraphics[scale=0.35, angle=-90, trim = 1.0cm 1.0cm 1.0cm 1.0cm, clip=true,]{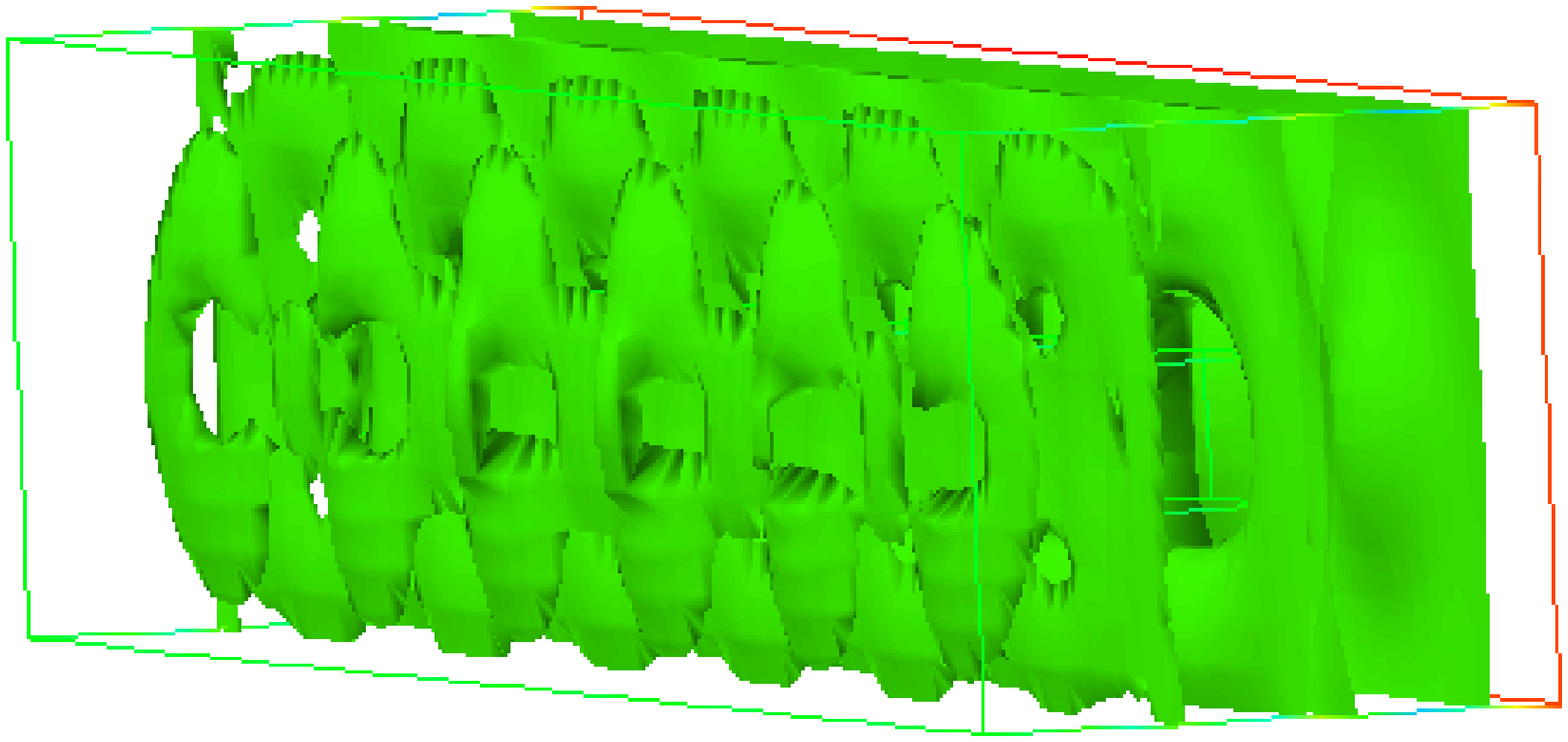}} \\
 e) $t=2.1$ & f) $t=2.1$  
 \end{tabular}
 \end{center}
 \caption{ Test 1. Isosurfaces of the simulated  FEM/FDM 
    solution of the problem (\ref{model1}) with initial conditions (\ref{initcond}) in $\Omega_{FEM}$  at different
    times. On a), c), d) we present transmitted data and on b), d), f) - backscattered data.}
 \label{fig:forward}
 \end{figure}

 \begin{figure}[tbp]
 \begin{center}
 \begin{tabular}{ccc}
 {\includegraphics[scale=0.5,  clip=]{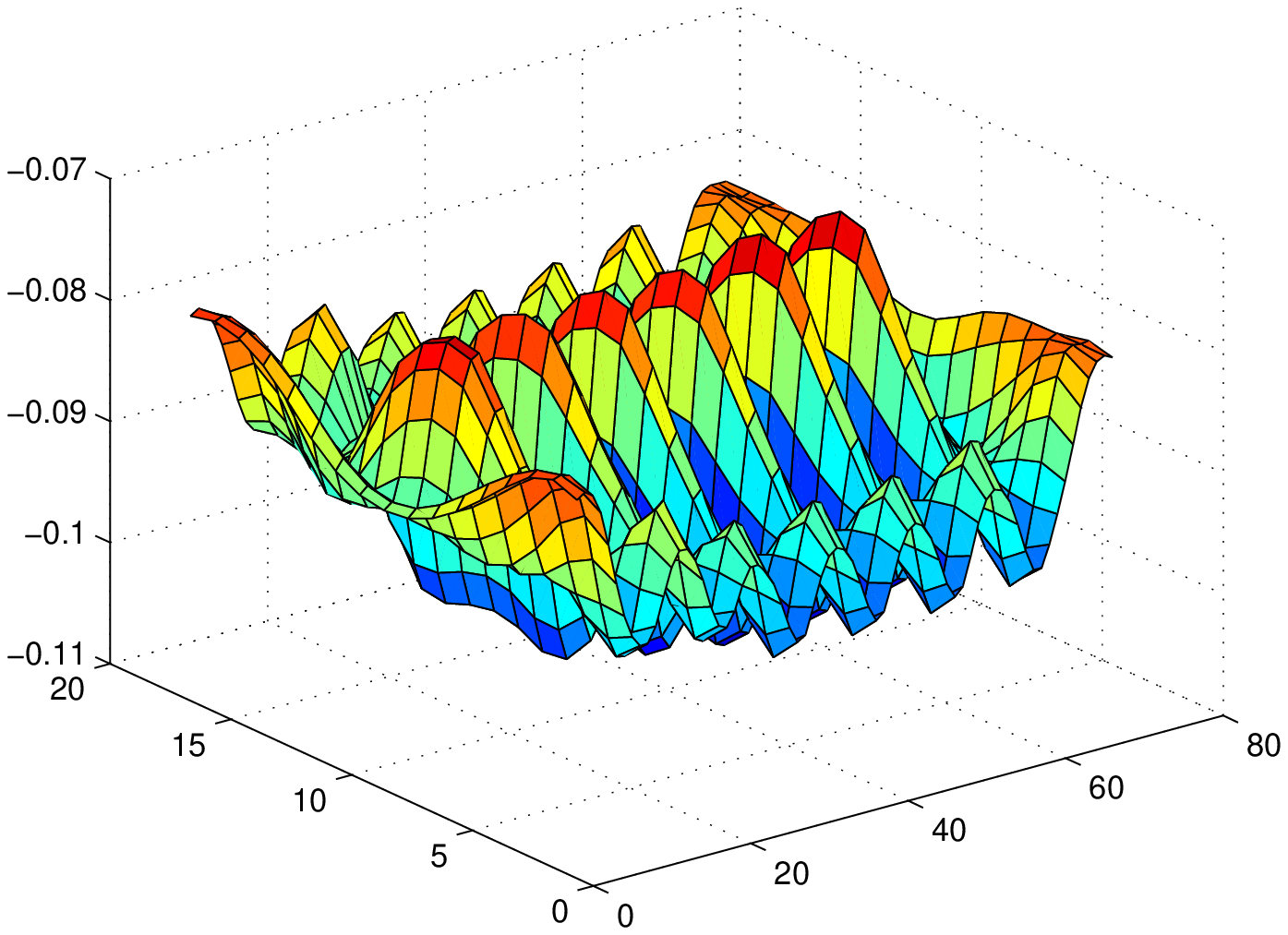}} &
 {\includegraphics[scale=0.5,  clip=]{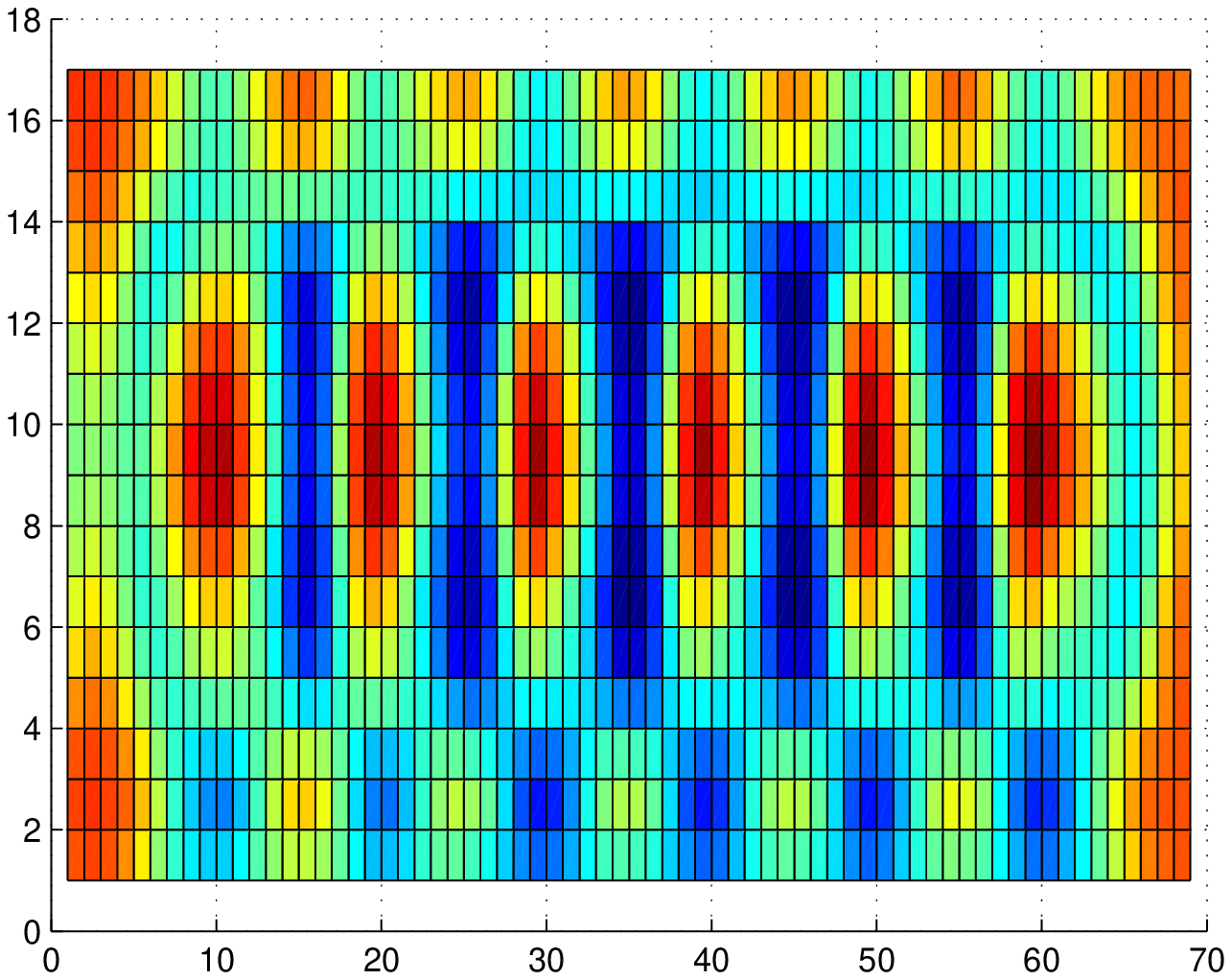}} \\
 a) Model 1 & b) $x_1 x_2$ view \\
 {\includegraphics[scale=0.5,  clip=]{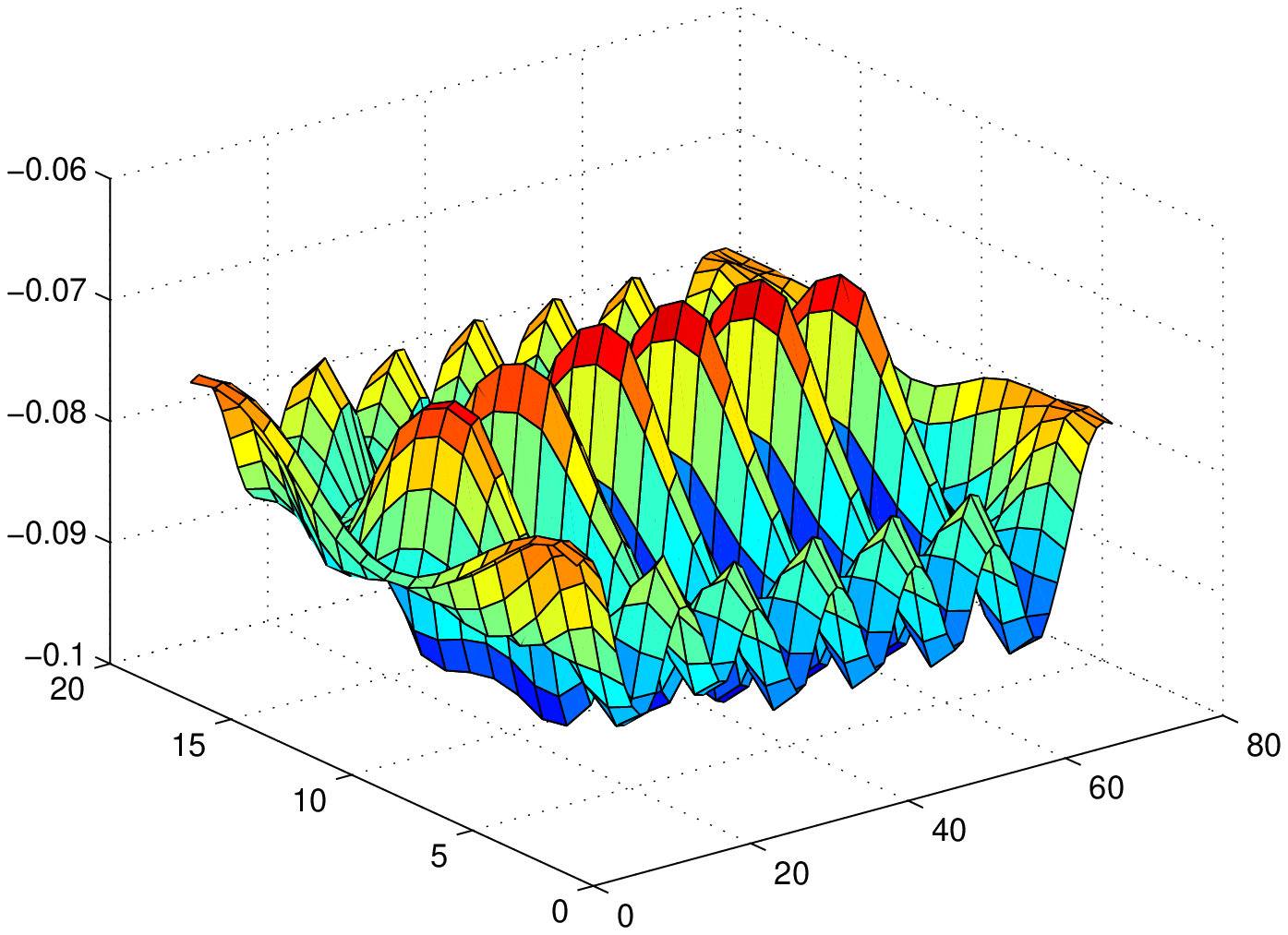}} &
 {\includegraphics[scale=0.5,  clip=]{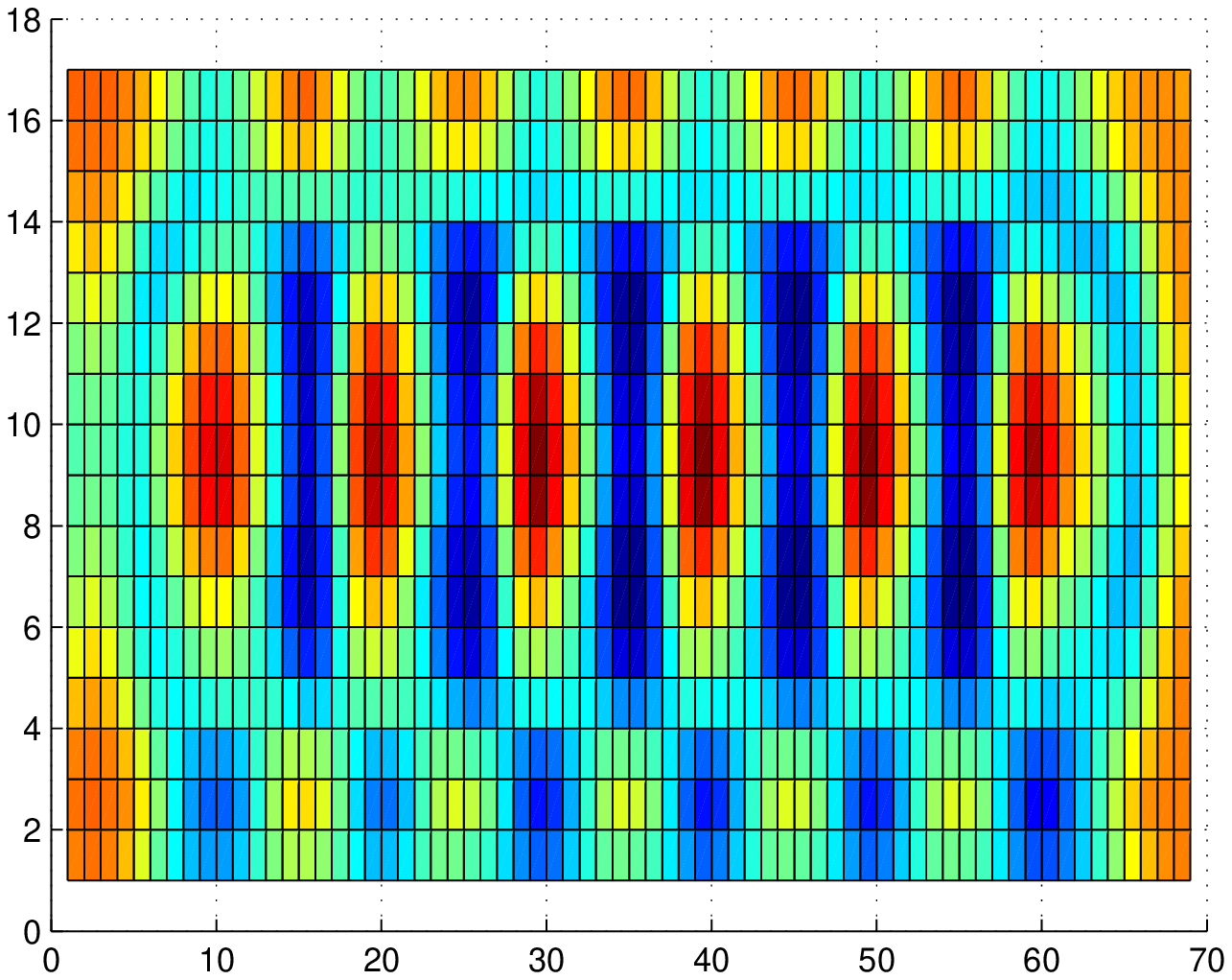}} \\
 c) Model 2 & d) $x_1 x_2$ view \\
 {\includegraphics[scale=0.5,  clip=]{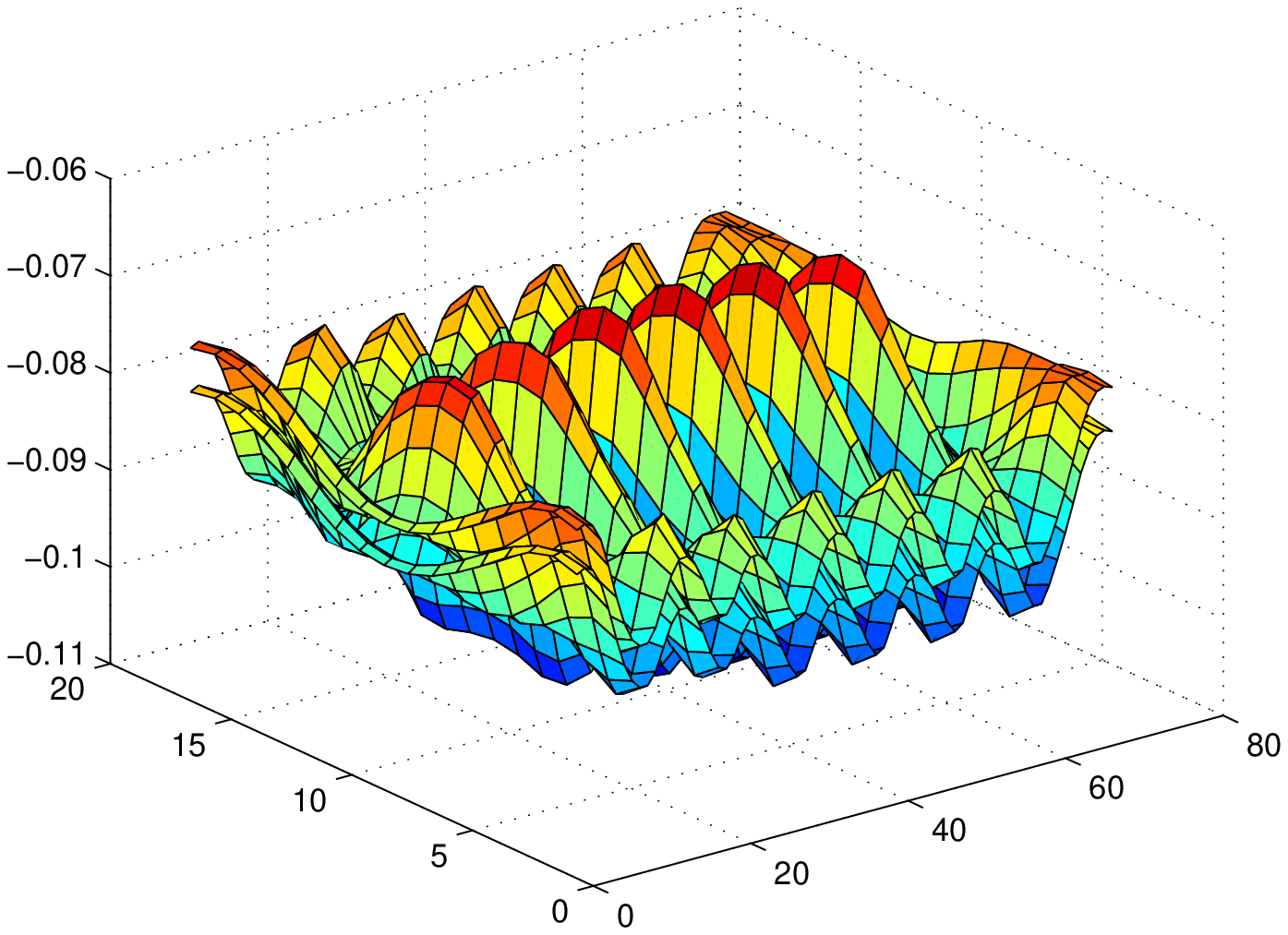}} &
 {\includegraphics[scale=0.5,  clip=]{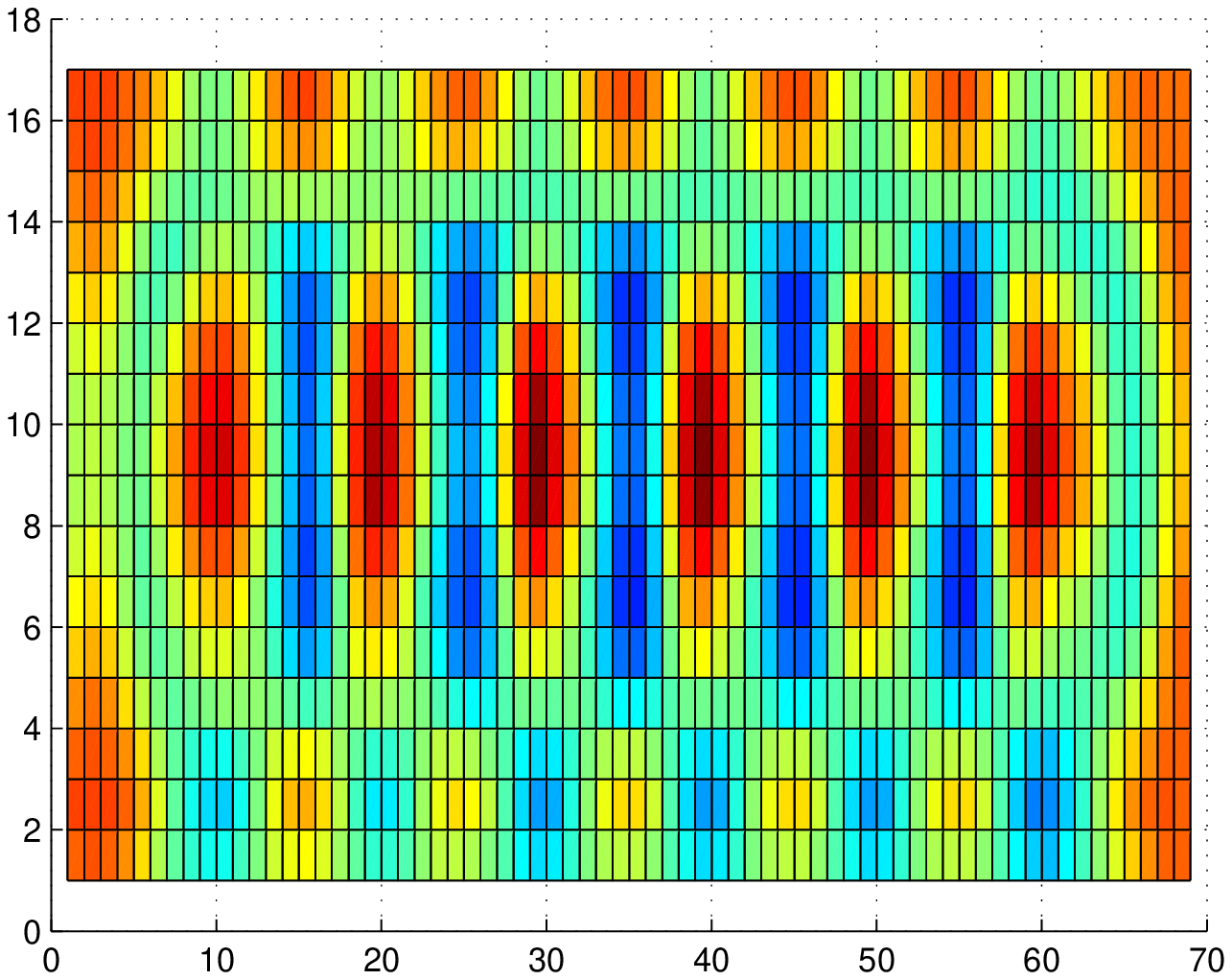}}
\\
 e) Models 1 and 2 & f) $x_1 x_2$ view \\
 \end{tabular}
 \end{center}
 \caption{ Test 1. Behavior of noisy backscattered data at time $t=1.8$ in both
   mathematical models of section \ref{sec:model}. Figure e) presents comparison of backscattered
   data in both models. }
 \label{fig:data}
 \end{figure}

 To generate backscattered data  at the observation points at $S_{T}$ in model problem 1,
 we solve the forward problem (\ref{model1}), with function $p(t)$
 given by (\ref{f}) in the time interval $t=[0,3.0]$ with the exact
 values of the parameters $a(x)=4.0$ inside
 scatterers of figure \ref{fig:fig1}, and $a(x)=1.0$
 everywhere else in $\Omega$.
We initialized initial conditions at backscattered side $\partial_1 \Omega$ as
\begin{equation}\label{initcond}
\begin{split}
u(x,0) &= \exp^{-(x_1^2 + x_2^2 + x_3^3)}  \cdot \cos  t|_{t=0} = \exp^{-(x_1^2 + x_2^2 + x_3^3)}  , \\
\frac{ \partial u}{\partial t} (x,0) &= -\exp^{-(x_1^2 + x_2^2 + x_3^3)} \cdot  \sin t|_{t=0} \equiv 0.
\end{split}
\end{equation}
  Figure \ref{fig:forward}-a) presents  behavior of this initial condition.

 Figure \ref{fig:data} presents typical behavior of noisy
 backscattered data for scatterers of figure \ref{fig:fig1} in our two
 models of section \ref{sec:model}.  Using figure \ref{fig:data}-e) we observe that the
 difference in the amplitude of these two data sets is very small, and
 thus, we expect that the influence of the non-zero initial conditions 
 will not affect to the reconstructions too much.

 Figure \ref{fig:fig4}-a) presents behavior of the computed $L_2$ norms
 of differences $F^m = ||\delta u^m - \delta u^{m-1}||$ for $m > 0$,
 where $\delta u^m = ||u^m - \tilde{u}|| z_{\delta}$ for $\omega=40$
 in (\ref{f}) and noise level $\sigma=3\%$. Analyzing
 this figure  for model problem 1 we observe that we
 achieve convergence in the optimization algorithm at iteration $m=10$
 in the conjugate gradient method. 
 Figures \ref{fig:fig4b} presents typical behavior of the computed
 $L_2$ norms of differences $F^m = ||\delta u^m - \delta u^{m-1}||$ for
 $m > 0$, where $\delta u^m = ||u^m - \tilde{u}|| z_{\delta}$ for
 different values of $\omega$ in (\ref{f}) and noise level $\sigma=10\%$.

We  can see results of reconstruction for model problem 1 in tables 1,2.
The reconstructed images of the conductivity function
 for both noise levels and different $\omega$  are presented in figures
 \ref{fig:fig7}, \ref{fig:fig9}.  Figures \ref{fig:fig7}-b) and \ref{fig:fig9}-c)  show best results of
 reconstruction which we have obtained  for $\omega=40$. We observe
 that for the noise $\sigma=3\%$ we get correct locations of scatterers and values of
 reconstructed parameter $a(x) \approx 4.73$ inside them compared with
 exact one $c=4.0$. For the noise $\sigma=10\%$ we get correct locations of scatterers and values of
 reconstructed parameter $a(x) \approx 4.35$ inside them.

Using these figures  we observe that the
location of all inclusions in $x_1 x_2$ directions is imaged very
well. However, from figure \ref{fig:fig10} follows
that the location in $x_3$ direction should still be improved.

\subsubsection{Test 2}

This test is similar to the previous one, only we solve \textbf{IP2}
in this case.  We start the optimization algorithm with guess values
of the parameters $a(x)=1.0$ at all points in $\Omega$. We use the
same as in (\ref{admpar}) set of admissible parameters $M_a$ and the
same regularization procedure as in test 1.

\begin{figure}[tbp]
 \begin{center}
 \begin{tabular}{ccc}
 {\includegraphics[scale=0.5,  clip=]{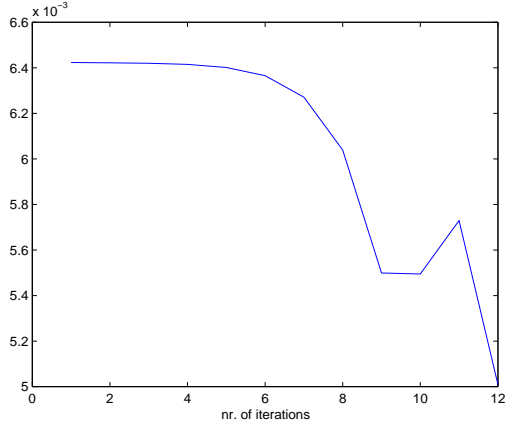}} &
 {\includegraphics[scale=0.5,  clip=]{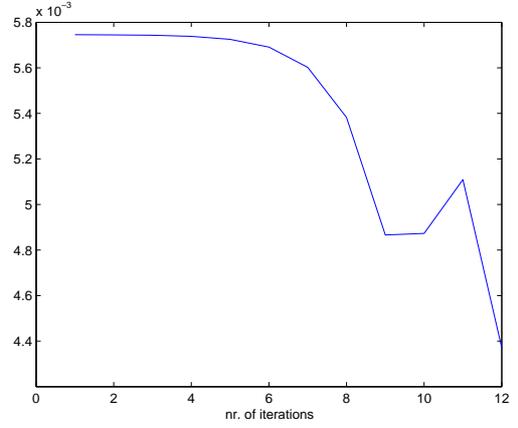}}
\\

a)  Model problem 1 & b)  Model problem 2
 \end{tabular}
 \end{center}
 \caption{ 
 Differences $F^m = ||\delta u^m - \delta u^{m-1}||$ 
   for $\omega=40$ in (\ref{f}).  Noise in backscattered data is 3\%.}
 \label{fig:fig4}
 \end{figure}

\begin{figure}[tbp]
 \begin{center}
 \begin{tabular}{cc}
{\includegraphics[scale=0.5,  clip=]{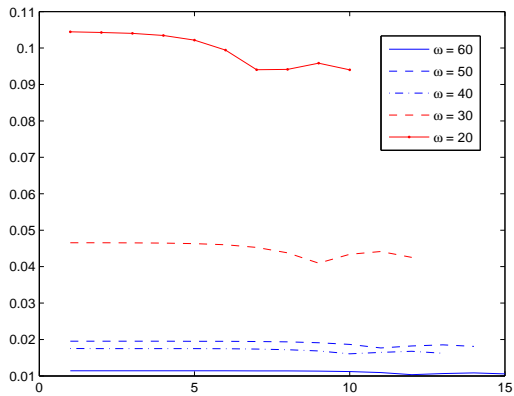}} &
 {\includegraphics[scale=0.5,  clip=]{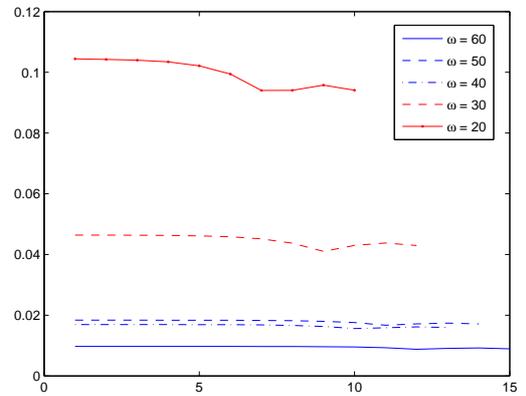}}\\
a) Model problem 1  & b) Model problem 2
 \end{tabular}
 \end{center}
 \caption{Behavior of differences $F^m = ||\delta u^m - \delta
   u^{m-1}||$ 
  for Model Problem 1 (left figures) and  for Model Problem 2 (right figures).  Noise in  backscattered data  is 10\%.}
 \label{fig:fig4b}
 \end{figure}

\begin{figure}[tbp]
 \begin{center}
 \begin{tabular}{cc}
 {\includegraphics[scale=0.35, angle=-90, trim = 6.0cm 1.0cm 7.0cm 6.0cm, clip=true,]{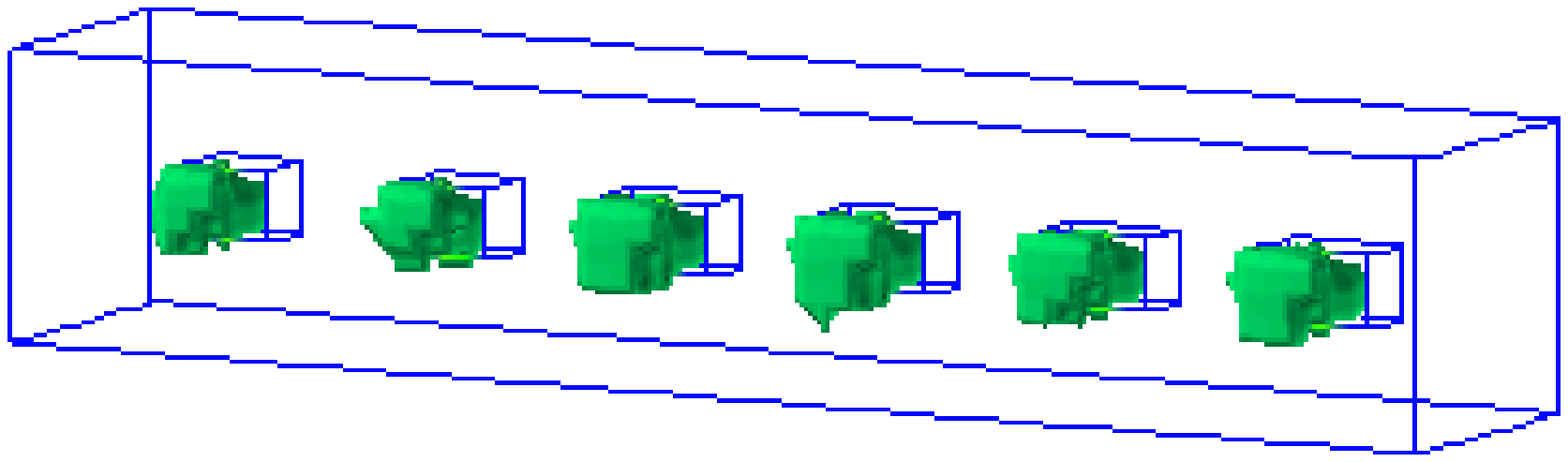}}  & 
 {\includegraphics[scale=0.35, angle=-90, trim = 6.0cm 1.0cm 7.0cm 6.0cm, clip=true,]{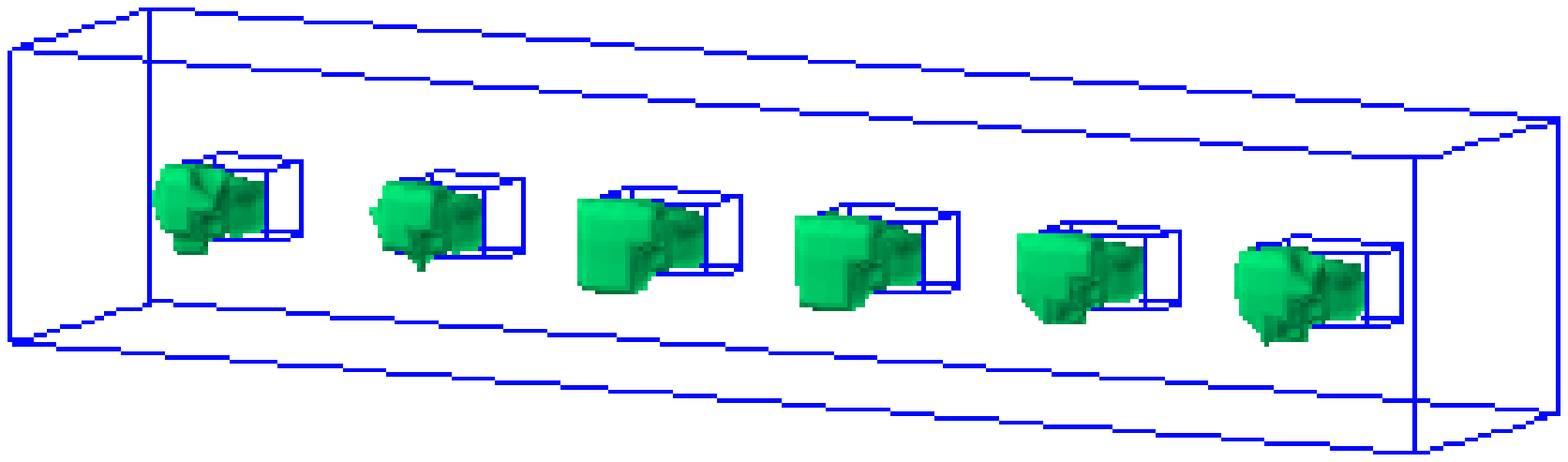}}
\\
a) $\omega=30,  \max\limits_{\Omega_{FEM} } a(x) = 5$ & b) $\omega=40,  \max\limits_{\Omega_{FEM} } a(x) = 4.73$  \\
{\includegraphics[scale=0.35, angle=-90, trim = 6.0cm 1.0cm 7.0cm 6.0cm, clip=true]{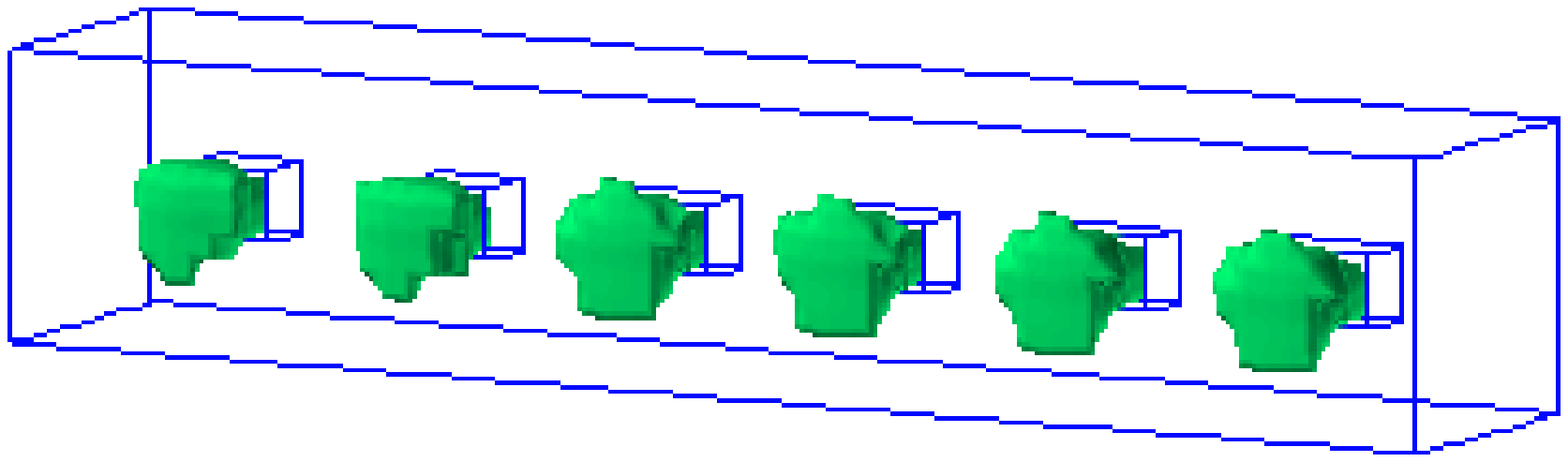}}  &
{\includegraphics[scale=0.35, angle=-90,  trim = 6.0cm 1.0cm 7.0cm 6.0cm, clip=true]{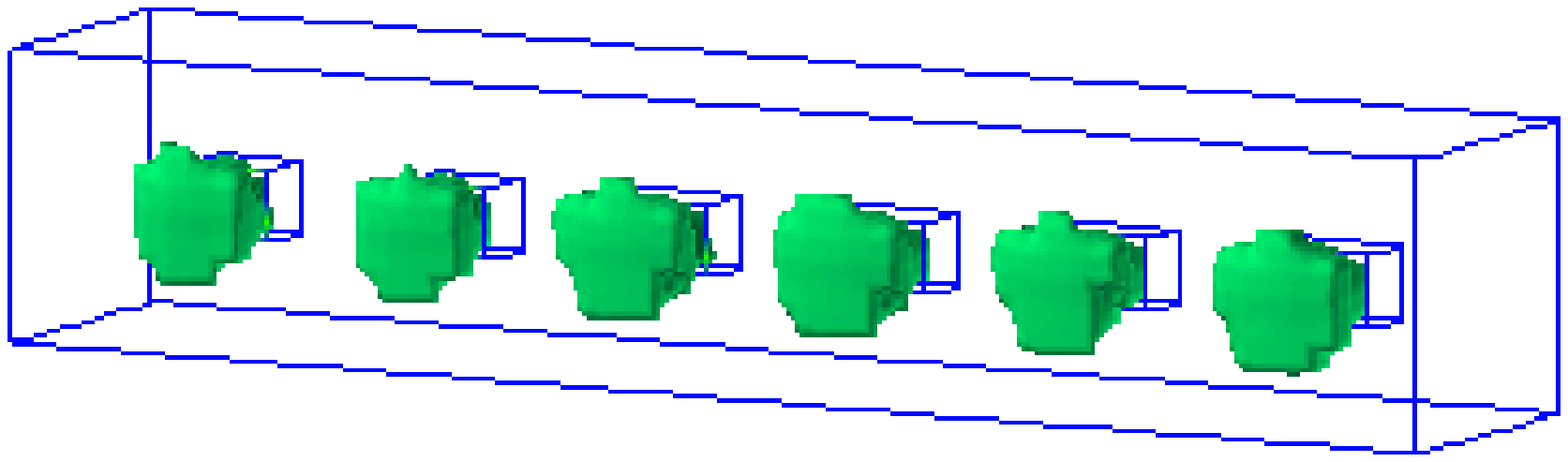}} 
\\
c) $\omega=50,  \max\limits_{\Omega_{FEM} } a(x) = 5.0$  & d) $\omega=60,  \max\limits_{\Omega_{FEM} } a(x) = 5.0$ 
 \end{tabular}
 \end{center}
 \caption{Test 1. Computed images of reconstructed functions $a(x)$ in model
   problem 1. We present functions $\tilde{a}$ for different $\omega$  in (\ref{f}) and noise level
   $\sigma=3\%$.}
 \label{fig:fig7}
 \end{figure}

Figure \ref{fig:fig4}-b) presents behavior of the computed $L_2$ norms
of differences $F^m = ||\delta u^m - \delta u^{m-1}||$, $m > 0$, for
$\omega=40$ in (\ref{f}) and noise level $\sigma=3\%$. Analyzing this
figure for model problem 2 we observe that we achieve convergence in
the optimization algorithm at iteration $m=10$ in the conjugate
gradient method.  Figure \ref{fig:fig4b}-b) presents typical behavior
of the computed $L_2$ norms of differences $F^m = ||\delta u^m -
\delta u^{m-1}||$ for $m > 0$, where $\delta u^m = ||u^m - \tilde{u}||
z_{\delta}$ for different values of $\omega$ in (\ref{f}) and noise
level $\sigma=10\%$ in this test.

\begin{figure}[tbp]
 \begin{center}
 \begin{tabular}{cc}
 {\includegraphics[scale=0.4, angle=-90, trim = 6.0cm 1.0cm 7.0cm 6.0cm, clip=true,]{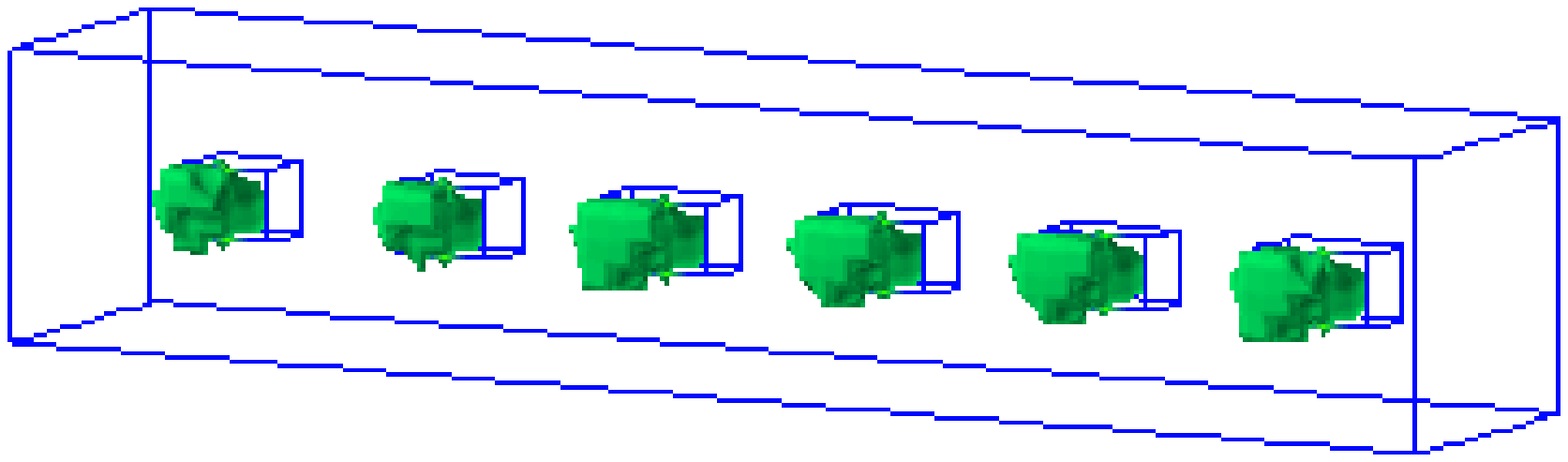}}  &
 {\includegraphics[scale=0.4, angle=-90, trim = 6.0cm 1.0cm 7.0cm 6.0cm, clip=true,]{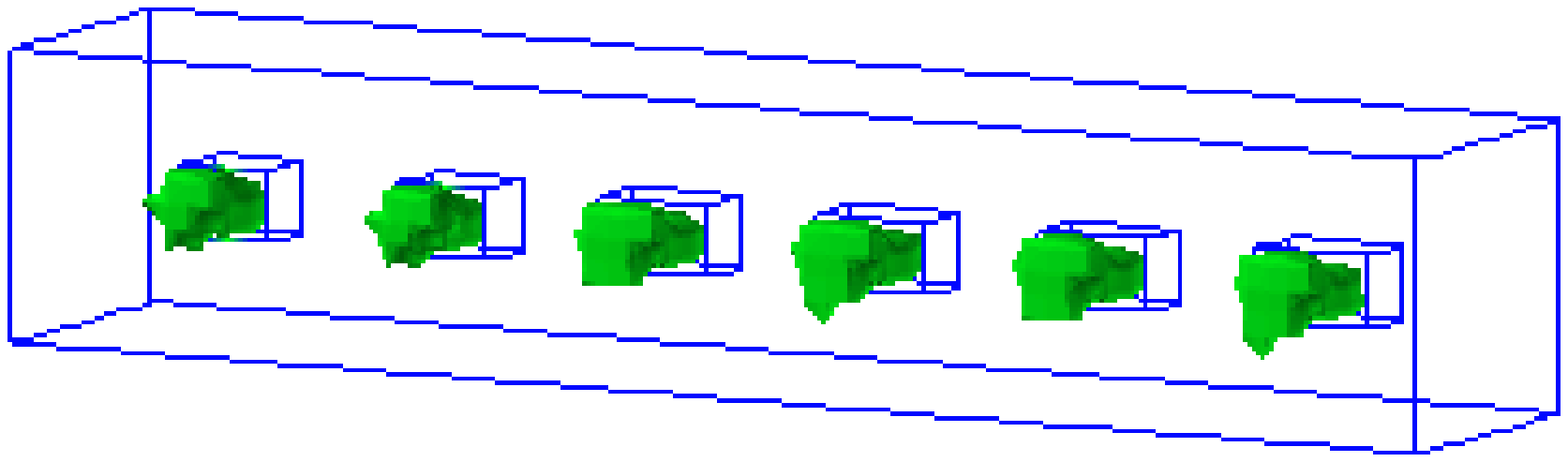}} 
 \\
a) $\omega=30,  \max\limits_{\Omega_{FEM} } a(x) = 4.86$  & b) $\omega=40,  \max\limits_{\Omega_{FEM} } a(x) = 4.19 $ \\
{\includegraphics[scale=0.4, angle=-90, trim = 6.0cm 1.0cm 7.0cm 6.0cm, clip=true]{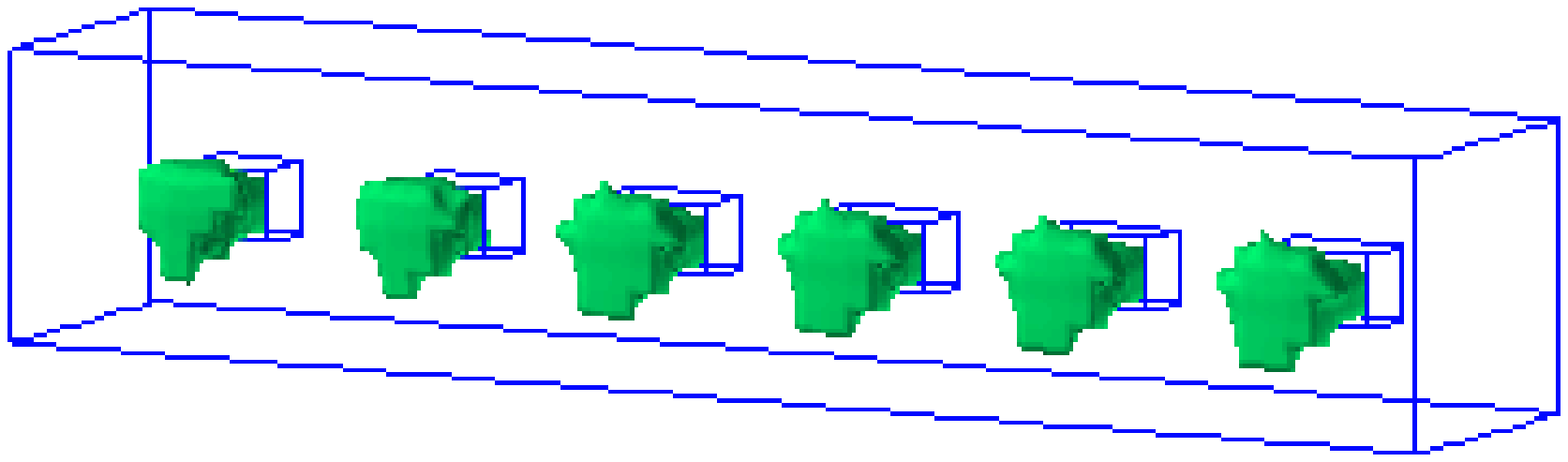}} &
{\includegraphics[scale=0.4, angle=-90,  trim = 6.0cm 1.0cm 7.0cm 6.0cm, clip=true]{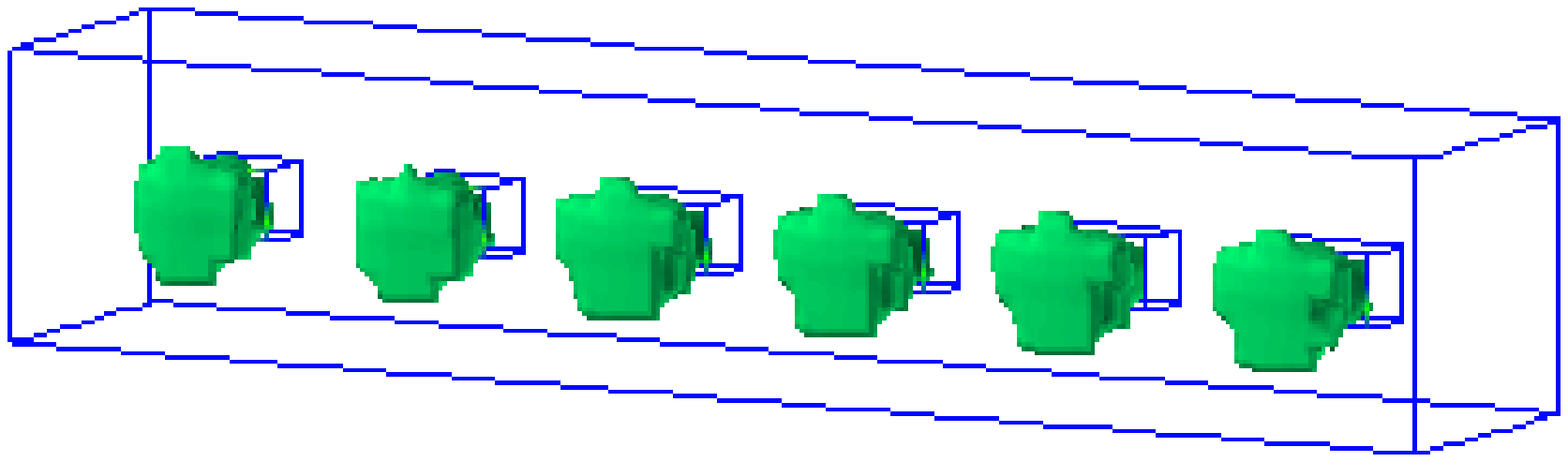}} 
\\
c) $\omega=50,  \max\limits_{\Omega_{FEM} } a(x) =5 $ &
d) $\omega=60,  \max\limits_{\Omega_{FEM} } a(x) =5 $ 
 \end{tabular}
 \end{center}
 \caption{ Test 2. Computed images of reconstructed functions $a(x)$ in model
   problem 2. We present functions $\tilde{a}$ for different $\omega$ in
   (\ref{f}) and noise level $\sigma=3\%$. Here we have initialized
   zero boundary conditions in the generation of backscattered data
   and in the optimization algorithm.}
 \label{fig:fig6}
 \end{figure}

\begin{figure}[tbp]
 \begin{center}
 \begin{tabular}{cc}
 {\includegraphics[scale=0.4, angle=-90, trim = 6.0cm 1.0cm 7.0cm 6.0cm, clip=true]{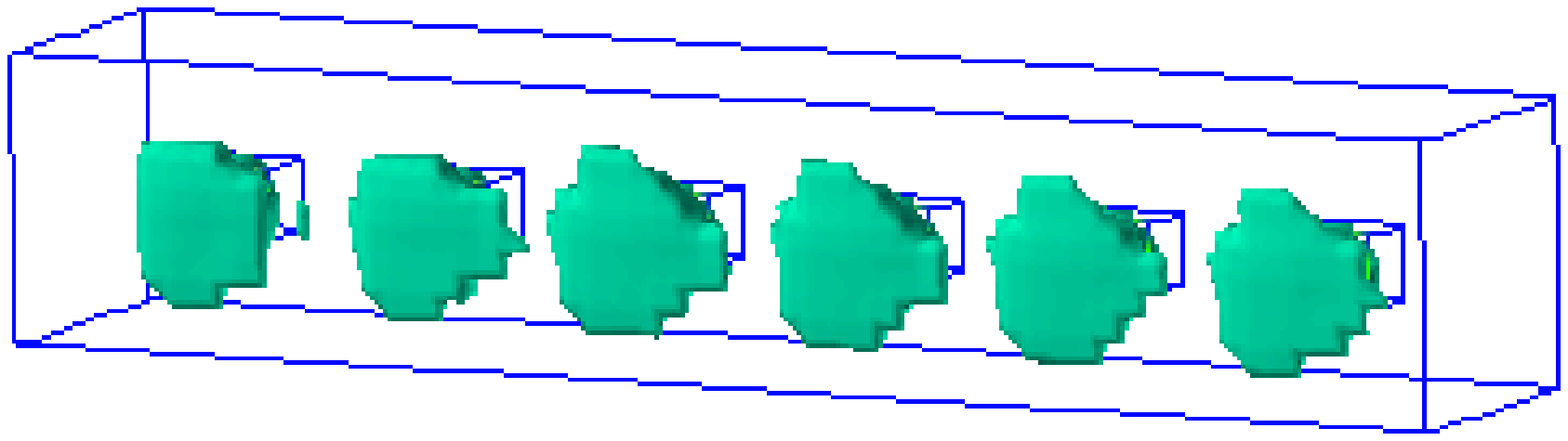}}  &
 {\includegraphics[scale=0.4, angle=-90, trim = 6.0cm 1.0cm 7.0cm 6.0cm, clip=true,]{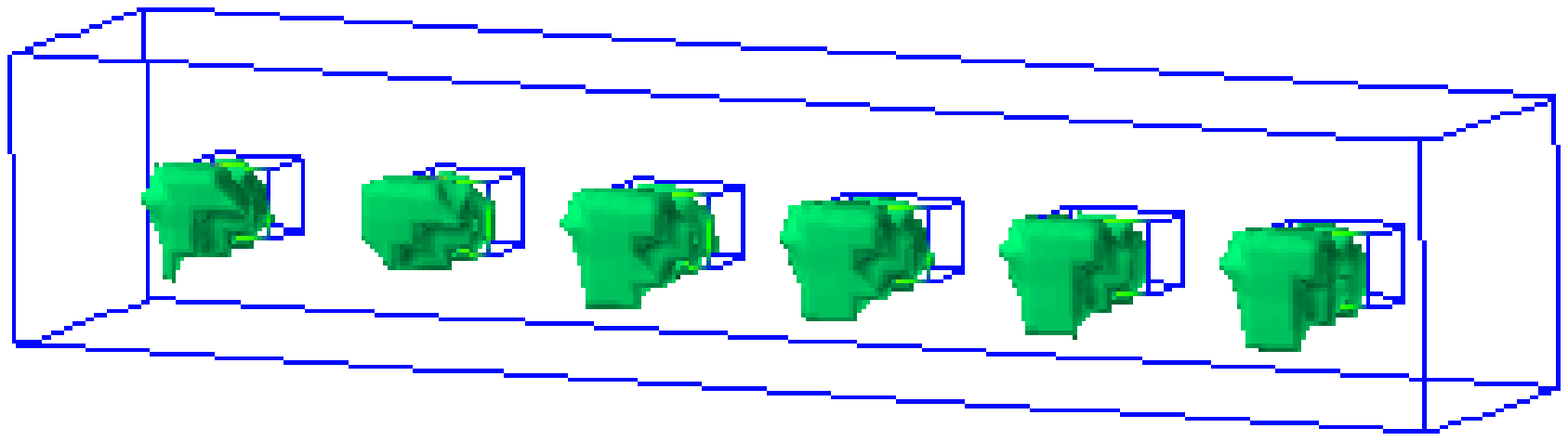}}
\\
a) $\omega=20,  \max\limits_{\Omega_{FEM} } a(x) = 3.29$  & b) $\omega=30,  \max\limits_{\Omega_{FEM} } a(x) = 4. 94$ \\
 {\includegraphics[scale=0.4, angle=-90, trim = 6.0cm 1.0cm 7.0cm 6.0cm, clip=true,]{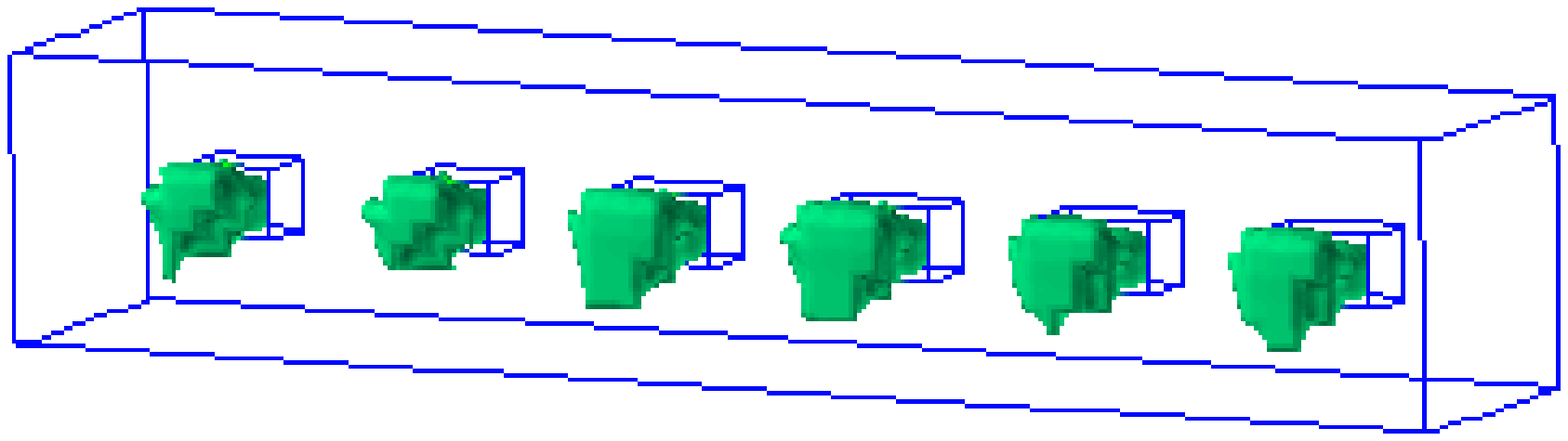}} &
{\includegraphics[scale=0.4, angle=-90, trim = 6.0cm 1.0cm 7.0cm 6.0cm, clip=true]{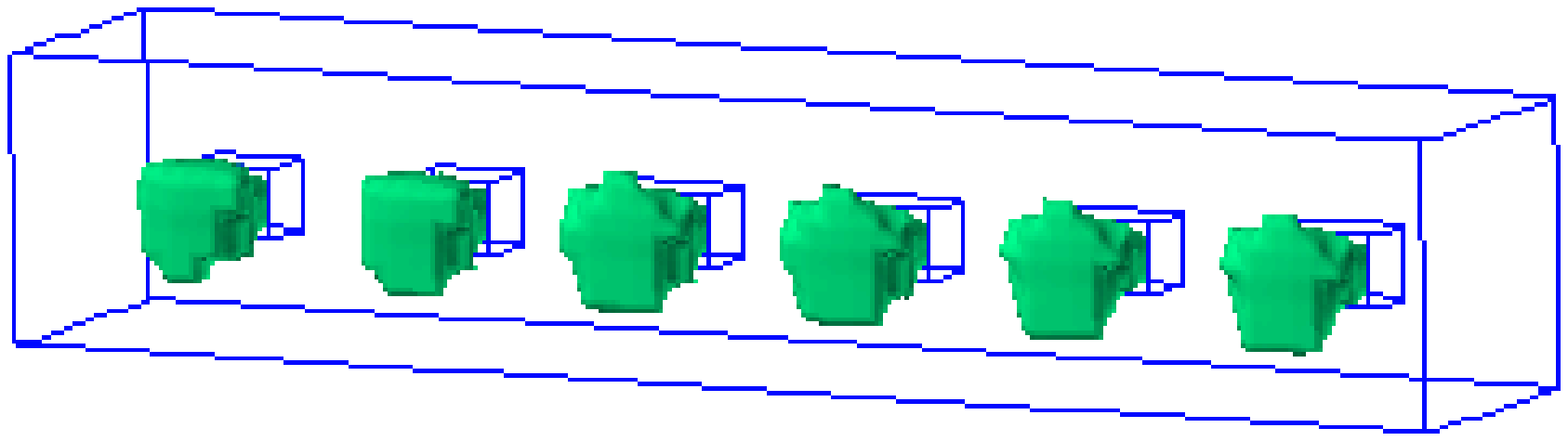}}
 \\
c) $\omega=40,  \max\limits_{\Omega_{FEM} } a(x) = 4.35$ &
d) $\omega=50,  \max\limits_{\Omega_{FEM} } a(x) = 4.4 $
 \end{tabular}
 \end{center}
 \caption{ Test 1. Computed images of reconstructed functions $a(x)$ in model
   problem 1. We present functions $\tilde{a}$ for different $\omega$ in
   (\ref{f}) and noise level $\sigma=10\%$.}
 \label{fig:fig9}
 \end{figure}

\begin{figure}[tbp]
 \begin{center}
 \begin{tabular}{cc}
 {\includegraphics[scale=0.35, angle=-90, trim = 6.0cm 1.0cm 7.0cm 6.0cm, clip=true,]{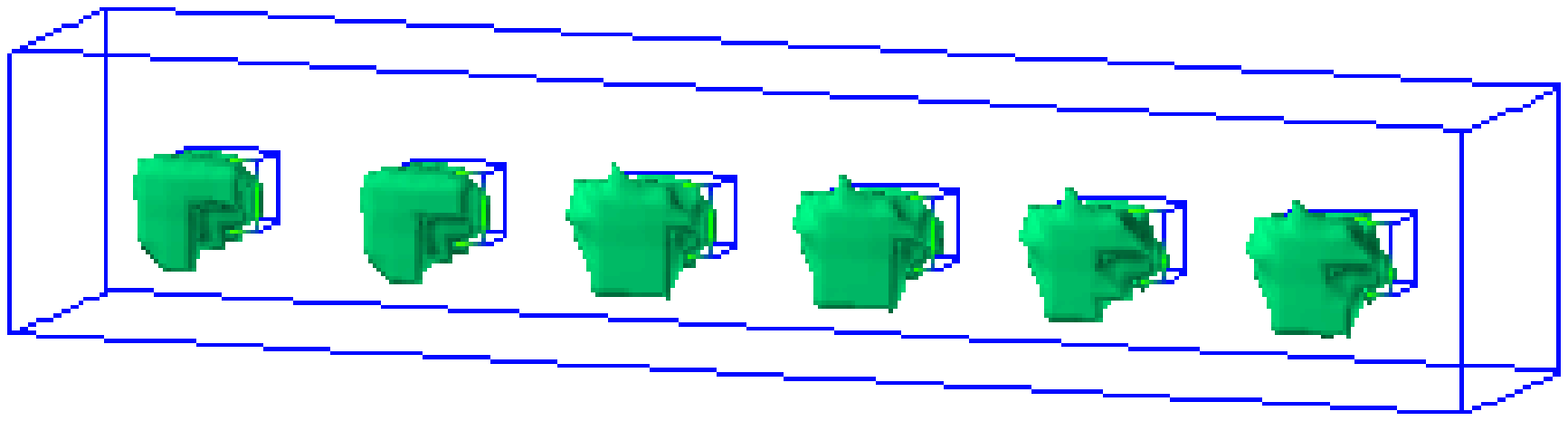}}  &
 {\includegraphics[scale=0.35, angle=-90, trim = 6.0cm 1.0cm 7.0cm 6.0cm, clip=true,]{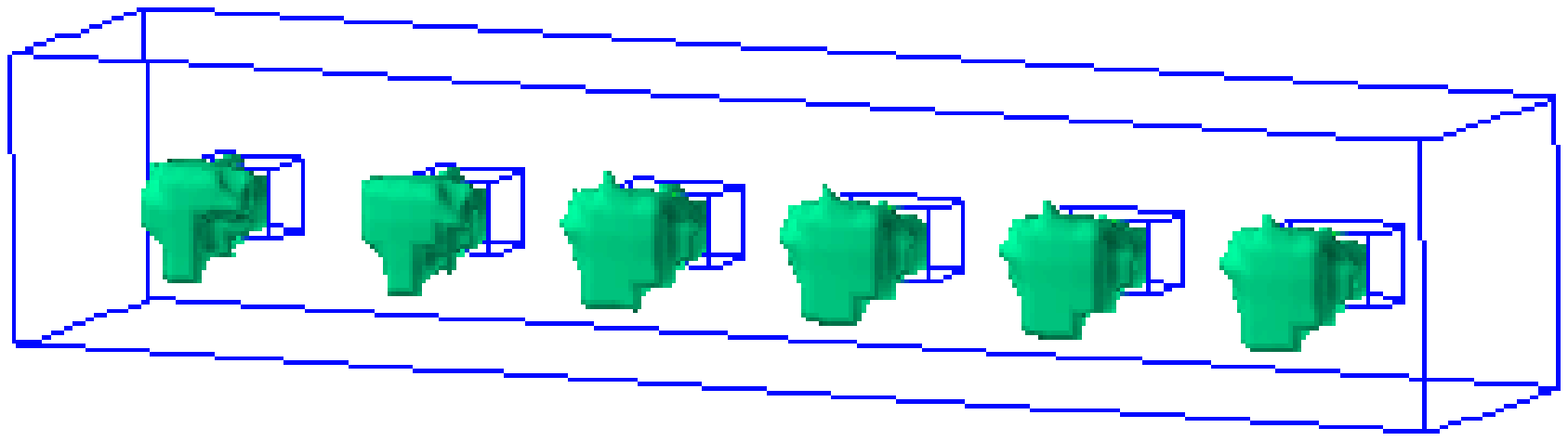}} 
\\
a) $\omega=30,  \max\limits_{\Omega_{FEM} } a(x) = 4.57$  & b) $\omega=40,  \max\limits_{\Omega_{FEM} } a(x) = 3.96 $ \\
{\includegraphics[scale=0.35, angle=-90, trim = 6.0cm 1.0cm 7.0cm 6.0cm, clip=true]{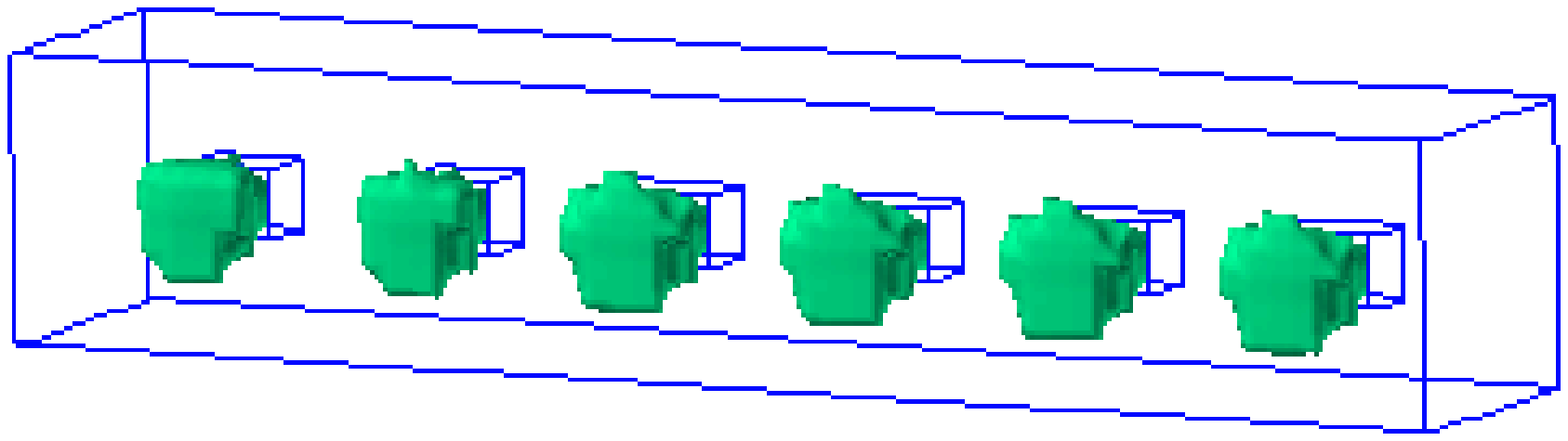}}  &
{\includegraphics[scale=0.35, angle=-90,  trim = 6.0cm 1.0cm 7.0cm 6.0cm, clip=true]{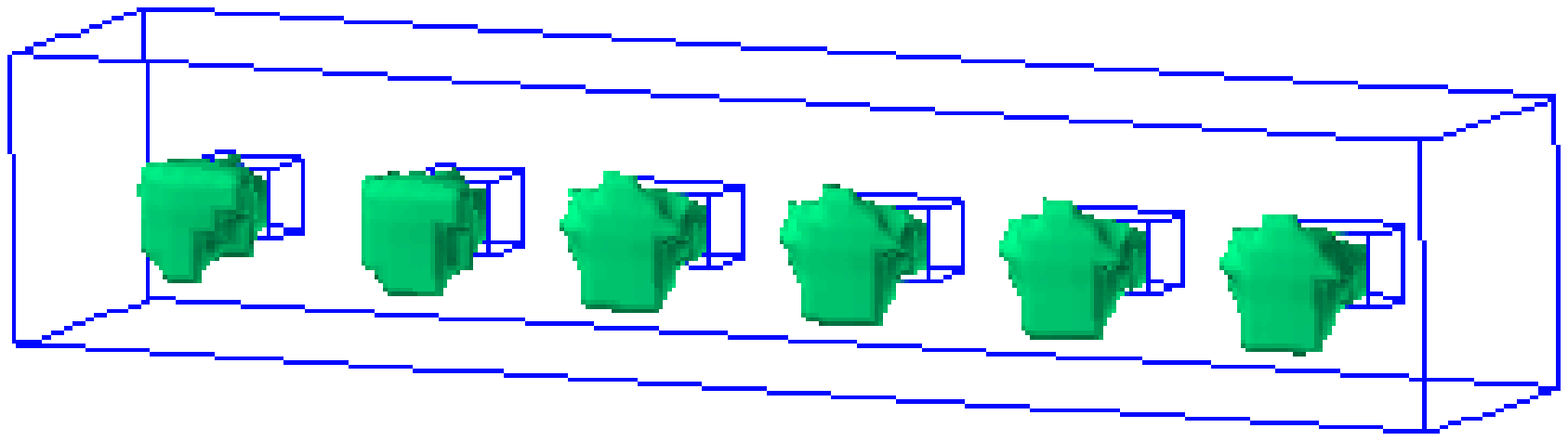}} 
\\
c) $\omega=50,  \max\limits_{\Omega_{FEM} } a(x) = 4.18 $  & d) $\omega=60,  \max\limits_{\Omega_{FEM} } a(x) = 4.46$  
 \end{tabular}
 \end{center}
 \caption{Test 2. Computed images of reconstructed functions $a(x)$ in model
   problem 2. We present functions $\tilde{a}$ for different $\omega$ in
   (\ref{f}) and noise level $\sigma=10\%$.}
 \label{fig:fig8}
 \end{figure}

Results of reconstruction of $a(x)$ for model problem 2 are presented
in tables 1,2.  Figures \ref{fig:fig6}, \ref{fig:fig8} show results of
reconstruction for both noise levels and different $\omega$.  Using
tables 1,2 we observe that best results of reconstruction are obtained
for $\omega=40$.  Figures \ref{fig:fig6}-b) and \ref{fig:fig8}-b) show
these results. Using figure \ref{fig:fig6}-b) we observe that for the
noise $\sigma=3\%$ we get correct locations of scatterers and values
of reconstructed parameter $a(x) \approx 4.19$ inside them compared
with exact one $a=4.0$ for $\omega=40$.  From figure \ref{fig:fig8}-b)
we see that for the noise $\sigma=10\%$ we get correct locations of
scatterers and values of reconstructed parameter $a(x) \approx 3.96$
inside them for $\omega=40$.

\begin{figure}[tbp]
 \begin{center}
 \begin{tabular}{cc}
 {\includegraphics[scale=0.3, angle=-90, trim = 1.0cm 1.0cm 1.0cm 1.0cm, clip=true,]{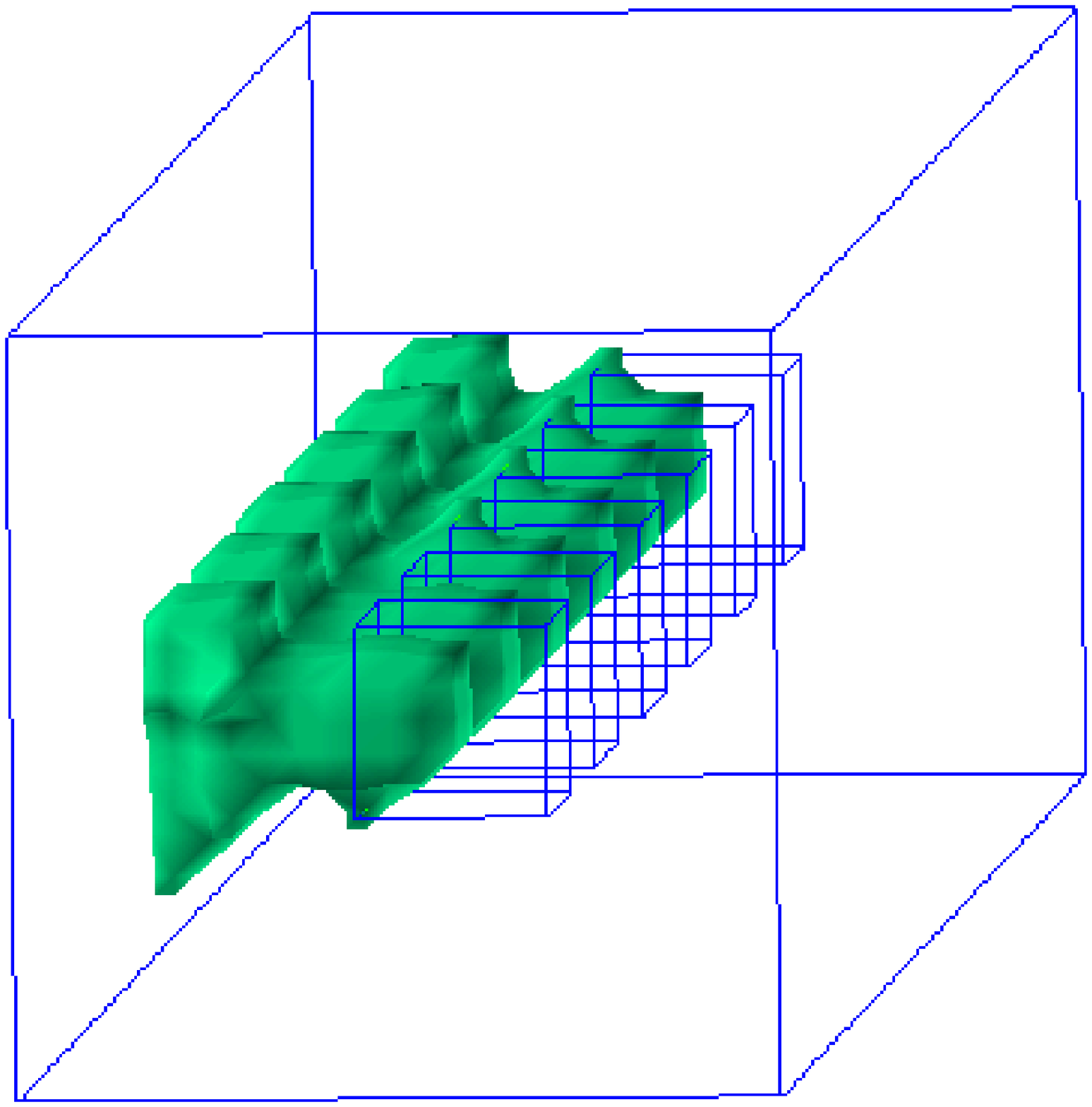}}  &
 {\includegraphics[scale=0.3, angle=-90, trim = 1.0cm 1.0cm 1.0cm 1.0cm, clip=true,]{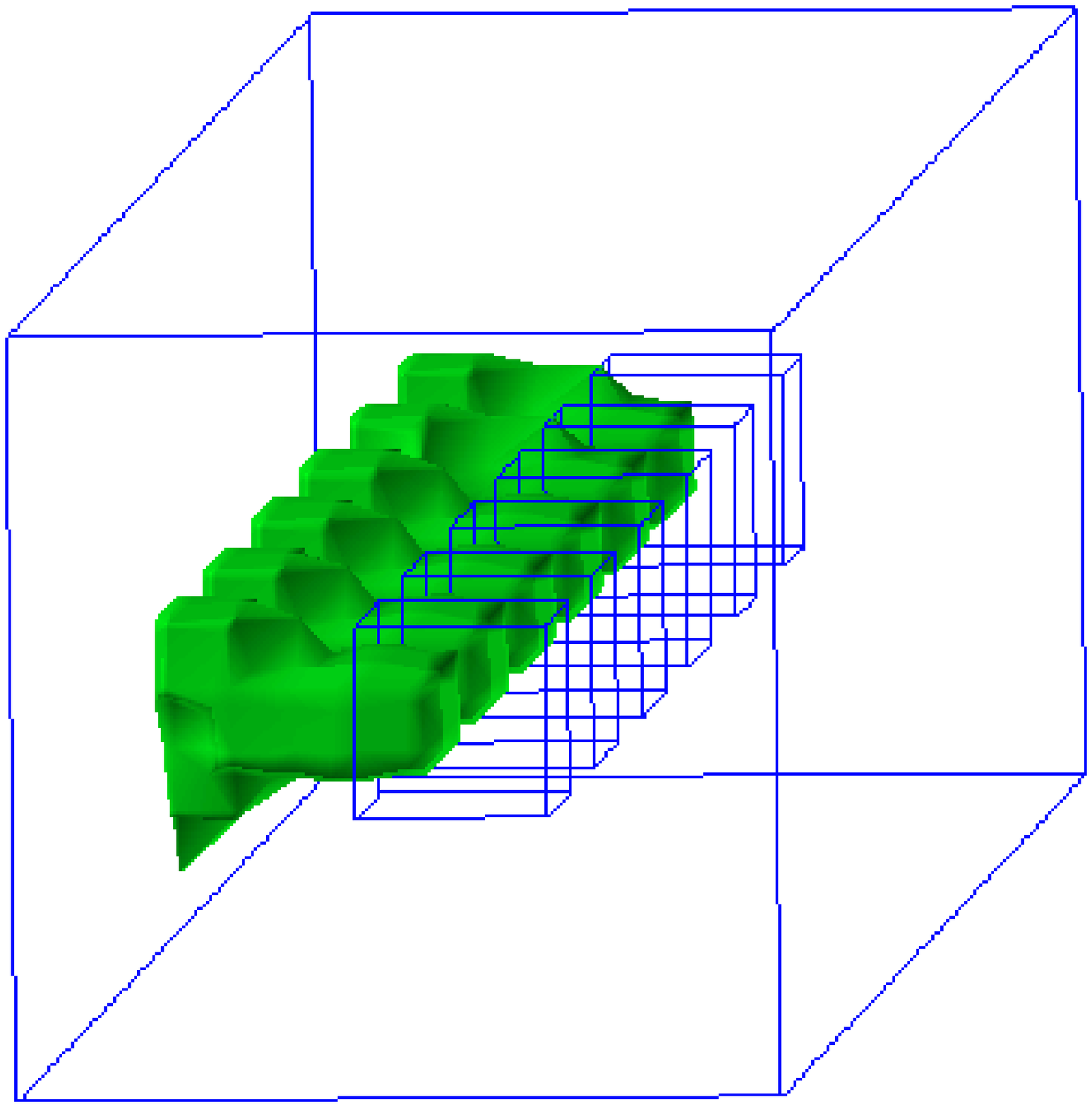}} 
\\
e) Test 1: $\omega=40, \sigma=10\%$  & f) Test 2: $\omega=40, \sigma =3\%$ \\
 \end{tabular}
 \end{center}
 \caption{ Computed images of reconstructed functions $a(x)$ in both
   model problems.  We present functions $\tilde{a}$ for 
   $\omega=40$ in (\ref{f}). We observe that
   reconstruction in $x_3$ direction should be improved.}
 \label{fig:fig10}
 \end{figure}

 Thus, we again conclude that the location of all
 inclusions in $x_1 x_2$ directions is imaged very well, but from
 figure \ref{fig:fig10} follows that the location in $x_3$
 direction should still be improved.

\section{Discussion and Conclusion}

In this work we present domain decomposition FEM/FDM method which is
applied for reconstruction of the conductivity function in the
hyperbolic equation in three dimensions. We have formulated inverse
problems and presented Lagrangian approach to solve these problems.
Explicit schemes for the solution of forward and adjoint problems in
the domain decomposition approach are also derived. We have
formulated different domain decomposition algorithms: the 
algorithm 1 describes the overlapping procedure between finite element
and finite difference domains, the algorithm 2 presents solution
of the forward and adjoint problems using the domain decomposition
FEM/FDM methods, and the  algorithm 3 describes the conjugate
gradient algorithm for reconstruction of the conductivity function
with usage of algorithms 1,2.

In our numerical tests we have obtained stable and good reconstruction
of the conductivity function $a(x)$ in the range of frequencies
$\omega \in [20,60]$ . Using tables 1,2 we can conclude that the best
reconstruction results are obtained in model problem 2 for $\omega=40$.  We can also
conclude that in all tests we could reconstruct size on $x_1 x_2$ -directions
for $a$, however, size in $x_3$ direction should be still
improved. Similarly with \cite{B, HybEf, BJ} we plan to apply an
adaptive finite element method in order to get better shapes and sizes
of the inclusions in all directions.

\section*{Acknowledgements}

This research is supported by the
Swedish Research Council (VR).
 The computations were performed on
resources at Chalmers Centre for Computational Science and Engineering
(C3SE) provided by the Swedish National Infrastructure for Computing
(SNIC).


\end{document}